\documentclass[leqno,12pt]{article}
\usepackage{amsmath,amssymb,rotating,a4wide,color}
\usepackage{wasysym}

\usepackage[all,cmtip]{xy}

\usepackage{verbatim}

\newdir^{ (}{{}*!/-3pt/\dir^{(}}
\newdir_{ (}{{}*!/-3pt/\dir_{(}}

\entrymodifiers={+!!<0pt,\fontdimen22\textfont2>}

\def\bigtimes{\mathop{\raise-2pt\hbox{\huge$\times$}}}

\newbox\circbulletbox
\setbox\circbulletbox=\hbox{\,{\large$\circledcirc$}\kern-8.7pt\raise .6pt\hbox{$\bullet$}\,}

\let\le\leqslant
\let\ge\geqslant

\def\mycirc{{\kern1pt\circ\kern2pt}}

\def\Im{\mathop{\rm Im}\nolimits}
\def\OI{{\overline{I}}}
\def\OU{{\overline{U}}}
\def\TI{\widetilde{I}}
\def\TJ{\widetilde{J}}

\def\Rad{\mathop{\rm Rad}\nolimits}

\def\Ker{\mathop{\rm Ker}\nolimits}

\def\proj{\mathop{\rm pr}\nolimits}

\let\phi\varphi
\let\theta\vartheta
\let\epsilon\varepsilon
\let\setminus\smallsetminus

\let\emptyset\varnothing

\newtheorem{Thm}{Theorem}[section]
\newtheorem{Prop}[Thm]{Proposition}
\newtheorem{Lem}[Thm]{Lemma}
\newtheorem{Claim}[Thm]{Claim}

\newtheorem{Def}[Thm]{Definition}

\newtheorem{Rem}[Thm]{Remark}
\newtheorem{Ex}[Thm]{Example}

\numberwithin{Thm}{section}
\def\UseTheoremCounterForNextEquation{\setcounter{equation}{\value{Thm}}\addtocounter{Thm}{1}}

\def\qed{{\hskip0pt\unskip\unskip\nobreak\hfil\penalty50
          \hskip1em\hbox{}\nobreak\hfil
           {$\square$}
          \parfillskip=0pt\finalhyphendemerits=0
          \par}\medskip}

\newenvironment{Proof}
               {\noindent{\bf Proof.}\ }
               {\qed}

               {\noindent{\bf Proof of #1.}\ }
               {\qed}

\newcommand{\BB}{{\mathbb{B}}}
\newcommand{\BC}{{\mathbb{C}}}

\newcommand{\BH}{{\mathbb{H}}}

\newcommand{\BP}{{\mathbb{P}}}

\newcommand{\BR}{{\mathbb{R}}}

\newcommand{\BZ}{{\mathbb{Z}}}

\newcommand{\Fm}{{\mathfrak{m}}}

\newcommand{\CC}{{\cal C}}

\newcommand{\UX}{{\underline{X}}}

\newcommand{\Uk}{{\underline{k}}}

\newcommand{\Ux}{{\underline{x}}}
\newcommand{\Uy}{{\underline{y}}}
\newcommand{\Uz}{{\underline{z}}}

\newbox\mybox
\def\arrover#1{\mathrel{
       \setbox\mybox=\hbox spread 1.4em
              {\hfil$\scriptstyle#1$\hfil}
       \vbox{\offinterlineskip\copy\mybox
             \hbox to\wd\mybox{\rightarrowfill}}}}

\def\larrover#1{\mathrel{
       \setbox\mybox=\hbox spread 1.4em
              {\hfil$\scriptstyle#1\vphantom{g}$\hfil}
       \vbox{\offinterlineskip\copy\mybox
             \hbox to\wd\mybox{\leftarrowfill}}}}

\def\ontoover#1{\mathrel{
       \setbox\mybox=\hbox spread 1.4em
              {\hfil$\scriptstyle#1\vphantom{g}$\hfil}
       \vbox{\offinterlineskip\copy\mybox
             \hbox to\wd\mybox{\rightarrowfill\hskip-2.8mm
                               $\rightarrow$}}}}
\def\leftontoover#1{\mathrel{
       \setbox\mybox=\hbox spread 1.4em
              {\hfil$\scriptstyle#1\vphantom{g}$\hfil}
       \vbox{\offinterlineskip\copy\mybox
             \hbox to\wd\mybox{$\leftarrow$\hskip-2.8mm
                               \leftarrowfill}}}}

\let\into\hookrightarrow

\begin{document}

\title{\strut
\vskip-80pt
The Strong Nullstellensatz\\
for Certain Normed Algebras}

\author{Richard Pink\\[12pt]
\small Department of Mathematics \\[-3pt]
\small ETH Z\"urich\\[-3pt]
\small 8092 Z\"urich\\[-3pt]
\small Switzerland \\[-3pt]
\small pink@math.ethz.ch\\[12pt]}



\maketitle

\begin{abstract}
Consider the polynomial ring in any finite number of variables over the complex numbers, endowed with the $\ell_1$-norm on the system of coefficients. Its completion is the Banach algebra of power series that converge absolutely on the closed polydisc. Whereas the strong Hilbert Nullstellensatz does not hold for Banach algebras in general, we show that it holds for ideals in the polynomial ring that are closed for the indicated norm. Thus the corresponding statement holds at least partially for the associated Banach algebra. We also describe the closure of an ideal in small cases.
\end{abstract}

{\renewcommand{\thefootnote}{}
}

\renewcommand{\baselinestretch}{0.75}\normalsize
\tableofcontents
\renewcommand{\baselinestretch}{1.0}\normalsize

\newpage


\section{Introduction}
\label{Intro}

{\bf Background:}
Fix a natural number $n$ and abbreviate $\BC[\UX] := \BC[X_1,\ldots,X_n]$. For any ideal $I\subset \BC[\UX]$ consider the zero set
\UseTheoremCounterForNextEquation
\begin{equation}\label{VIDef}
V(I)\ :=\ \bigl\{ \Uz\in\BC^n \bigm| \forall f\in I\colon f(\Uz)=0\bigr\}.
\end{equation}
Dually, for any subset $Z\subset\BC^n$ consider the vanishing ideal
\UseTheoremCounterForNextEquation
\begin{equation}\label{IZDef}
I(Z)\ :=\ \bigl\{ f\in\BC[\UX] \bigm| \forall \Uz\in Z\colon f(\Uz)=0\bigr\}.
\end{equation}
Furthermore, the radical of an ideal $I$ is the ideal
\UseTheoremCounterForNextEquation
\begin{equation}\label{RadDef}
\Rad(I)\ :=\ \bigl\{ f\in\BC[\UX] \bigm| \exists k\ge1\colon f^k\in I\bigr\}.
\end{equation}
\emph{Hilbert's Nullstellensatz} \cite{Hilbert1893} asserts that for any ideal $I$ we have $I(V(I))=\Rad(I)$.
The special case with $V(I)=\emptyset$ is often called the \emph{weak Nullstellensatz} and says that $V(I)=\emptyset$ if and only if $I=(1)$. 
One strategy of proof, due to Rabinowitsch \cite{Rabinowitsch1930}, is to deduce the strong Nullstellensatz from the weak one using localization (see, e.g., Harris \cite[Ch.\;5]{HarrisAlgGeom1992}). 

Both results are really properties of the finitely generated $\BC$-algebra $A:=\BC[\UX]/I$. Namely, the set $V(I)$ is in natural bijection with the set of all $\BC$-algebra homomorphisms $A\to\BC$, and the weak Nullstellensatz says that $V(I)=\emptyset$
if and only if $A=0$, while the strong Nullstellensatz says that any element of $A$ whose value at each point of $V(I)$ is zero is nilpotent.

\medskip
Now let $A$ be a commutative complex Banach algebra, and let $M(A)$ denote the set of all continuous $\BC$-algebra homomorphisms $A\to\BC$. A basic result (Naimark \cite[\S9.4, \S11.1]{Naimark1972}) 
states that there is a natural bijection from $M(A)$ to the set of all maximal ideals of $A$ defined by $\phi\mapsto\Ker(\phi)$; in particular, all maximal ideals of $A$ are closed. From this one deduces the \emph{weak Nullstellensatz for commutative Banach algebras}, namely that $M(A)=\emptyset$ if and only if $A=0$.

The analogue of the strong Nullstellensatz would be the statement that any element of~$A$, whose value at each point of $M(A)$ is zero, is nilpotent. This, however, is false in general. 
But the usual counterexamples are based on deliberately devious constructions (e.g. \cite{RadicalBanachAlgebras1983})
and do not arise naturally in algebraic geometry. One may therefore reasonably ask whether a version of the strong Nullstellensatz still holds for the commutative Banach algebras that one usually encounters. 
This article provides a partial answer to this question.

\medskip
Note that literature on Banach algebras often leaves out the adjective `weak' when referring to the weak Nullstellensatz. Usually the goal is either to give an elementary proof of the weak Nullstellensatz in a special situation, such as in 
von Renteln \cite{vonRenteln1981}, 
Mortini-von Renteln \cite{MortiniVonRenteln1989}, 
Bridges-Mines-Richman-Schuster \cite{BridgesMinesRichmanSchuster2003}, 
Mortini-Rupp \cite{MortiniRupp2012},
or to prove an analogue of the weak Nullstellensatz dealing only with an open part of the spectrum, as in Gelca \cite{Gelca1995}, \cite{Gelca1997}, or in connection with the corona problem as in Carleson \cite{Carleson1962}, Krantz-Li \cite{KrantzLi1996}, and many others. Namely, via the weak Nullstellensatz the corona statement is equivalent to saying that a certain open subset of the spectrum of a Banach algebra is dense, and so possesses no corona.

\medskip
{\bf Results of this article:}
We consider $\BC[\UX]$ as a normed $\BC$-algebra by setting
\UseTheoremCounterForNextEquation
\begin{equation}\label{NormDef}
\textstyle \bigl\Vert\sum_\Uk a_\Uk \UX^\Uk\bigr\Vert \ :=\ \sum_\Uk |a_\Uk| \ \in\ \BR^{\ge0}.
\end{equation}
The reason for this choice is the following universal property:

\begin{Prop}\label{UP}
\begin{enumerate}
\item[(a)] For any normed $\BC$-algebra $(A,\Vert\ \Vert)$ and any elements $a_j\in A$ satisfying ${\Vert a_j\Vert\le1}$, there exists a unique $\BC$-algebra homomorphism $\phi\colon \BC[\UX]\to A$ such that ${\phi(X_j)=a_j}$ for all $1\le j\le n$ and ${\Vert\phi(f)\Vert\le\Vert f\Vert}$ for all $f\in\BC[\UX]$.
\item[(b)] The norm (\ref{NormDef}) is the only norm on $\BC[\UX]$ which has the universal property (a) and satisfies $\Vert X_j\Vert=1$ for all $1\le j\le n$.
\end{enumerate}
\end{Prop}

\begin{Proof}
Combine the universal property of polynomial rings with the defining properties of norms on algebras.
\end{Proof}

Let $\BB:=\{z\in\BC:|z|\le1\}$ denote the closed unit disc and $\partial\BB = \{z\in\BC:|z|=1\}$ its boundary, the unit circle. 
The completion of $\BC[\UX]$ with respect to the chosen norm is the algebra of power series that converge absolutely on the closed polydisc~$\BB^n$. In fact, this completion is a Banach algebra with spectrum naturally homeomorphic to~$\BB^n$. In order to stay closer to algebraic geometry we will, however, continue working with the normed algebra $\BC[\UX]$.

For any ideal $I\subset \BC[\UX]$ consider the restricted zero set
\UseTheoremCounterForNextEquation
\begin{equation}\label{MIDef}
M(I)\ :=\ \bigl\{ \Uz\in\BB^n \bigm| \forall f\in I\colon f(\Uz)=0\bigr\}
\ =\ V(I)\cap \BB^n.
\end{equation}
As usual, we denote the closure of $I$ by~$\OI$, which is again an ideal (Naimark \cite[\S8.1]{Naimark1972}). Since the evaluation map $\BC[\UX]\to\BC$ at any point of $\BB^n$ is continuous with respect to the norm $\Vert\ \Vert$, we have $M(I)=M(\OI)$. The first main result of this article is the following analogue of the strong Nullstellensatz:

\begin{Thm}\label{NS1}
For any ideal $I\subset\BC[\UX]$ we have $I(M(I)) = \Rad(\OI)$.
\end{Thm}

This result can be interpreted as giving some information about the closure~$\OI$, but not all. It leads to the question whether one can describe $\OI$ precisely in a purely algebraic manner. To this we give the following partial answers.
For any point $\Uz=(z_1,\ldots,z_n)\in\BC^n$ consider the maximal ideal $\Fm_\Uz := (X_1-z_1,\ldots,{X_n-z_n})$ of $\BC[\UX]$.

\begin{Thm}\label{NS2}
For any $\Uz\in\BC^n$ and any $\Fm_\Uz$-primary ideal $I\subset\BC[\UX]$ we have
$$\OI\ =\ \
\biggl\{\!\begin{array}{ll}
I + (X_j-z_j)_{1\le j\le n,\, |z_j|=1}
 & \hbox{if $\Uz\in\BB^n$,} \\[3pt]
(1) & \hbox{if $\Uz\not\in\BB^n$.}
\end{array}$$
\end{Thm}

We can interpret this result heuristically as saying that only points $\Uz\in\BB^n$ contribute to the spectrum of $(\BC[\UX],\Vert\ \Vert)$ and that within this spectrum, infinitesimal deformations in the direction $X_j$ are possible to arbitrarily high degree if $z_j\in\BB^\circ$, but not at all if $z_j\in\partial\BB$.
The last statement is geometrically plausible in so far as any algebraic deformation of $z_j\in\partial\BB$ would have tangent space $\BC$ and would thus include an infinitesimal deformation in the direction out of~$\BB$, which should be impeded.

We also analyze the situation in dimension~$2$, the most interesting case being:

\begin{Thm}\label{NS3}
Consider any irreducible polynomial $f\in\BC[X_1,X_2]$ such that $M((f))$ is infinite and contained in $(\partial\BB)^2$. Then for any $k\ge1$ we have 
$$\overline{(f^k)}\ =\ (f).$$
\end{Thm}

The heuristic explanation for this result is the same as above. 
We also have a complete answer in dimension $\le2$ or in the case of finite support. Note that Theorem \ref{NS2} yields an explicit description of $\overline{\Fm_\Uz^r}$ for all $\Uz\in\BB^n$ and $r\ge1$. 

\begin{Thm}\label{NS4}
Consider any ideal $I\subset\BC[\UX]$. If $M(I)$ is finite or $n\le 2$, then 
$$\OI\ =\ \bigcap_{\Uz\in\BB^n} \bigcap_{r\ge1} \bigl(I+\overline{\Fm_\Uz^r}\,\bigr).$$
\end{Thm}

I am sorely tempted to conjecture that Theorem \ref{NS4} is true without any condition on $I$ or~$n$. However, at present I cannot exclude the possibility that the ambient dimension $n$ has some influence. 
For example consider the polynomial
$$f\ :=\ 1+w\cdot(X_1+X_2+X_3+X_1X_2+X_1X_3+X_2X_3)+X_1X_2X_3 \ \in \ \BC[X_1,X_2,X_3]$$
for a real number $0<w<1$. Its associated $M((f))$ is contained in $(\partial\BB)^3$ and Zariski dense in the hypersurface $V(f)$. Therefore $\overline{(f)}=(f)$. With the methods used to prove Theorem \ref{NS3}, I can show that $f^2\in\overline{(f^k)}$ for all $k\ge2$, but not the remaining step $f\in\overline{(f^2)}$.

\medskip
{\bf Methods:}
The proof of Theorem \ref{NS1} is contained in Sections \ref{Prelim} and \ref{NSProof}. Consider a system of generators $(f_1,\ldots,f_m)$ of~$I$ and a polynomial $g\in I(M(I))$. We first use the {\L}ojasiewicz inequality (see Bierstone-Milman \cite[Thm.\;6.4]{BierstoneMilman1988}) to write some power of $g$ as a linear combination $g^N=\sum_{j=1}^m g_jf_j$ with functions $g_j\in C^M(\BB^n)$ for sufficiently large~$M$. Then we improve this representation successively by explicitly solving the $\smash{\bar\partial}$-equation as in Krantz-Li \cite[\S2]{KrantzLi1996}. Hidden behind this is really the Koszul complex (compare Costea-Sawyer-Wick \cite{CosteaSawyerWick2011}). This construction ends with another such linear combination, where the $g_j$ are in addition holomorphic on $(\BB^\circ)^n$. If now $M>\frac{n}{2}$, those functions are represented by power series which converge absolutely on~$\BB^n$. Approximating these by polynomials finally shows that $g^N\in\OI$.

In Section \ref{FinSupp} we deal with the case of finite support and prove Theorem \ref{NS2} and the finite support part of Theorem \ref{NS4}, using relatively direct calculations.

In Sections \ref{PlaneCurves} and \ref{PlaneD} we analyze plane curves in some  detail and prove Theorem \ref{NS3} as part of Theorem \ref{NS3a}. 
The key point here is Proposition \ref{NS3aMainProp}, establishing something like an approximate identity for the ideal $(f)$ (compare Mortini-von Renteln \cite[page 223]{MortiniVonRenteln1989}).

Section \ref{Dim2} establishes the case $n=2$ of Theorem \ref{NS4} by combining the same kind of arguments as before with a certain amount of commutative algebra.

We end this article with some examples in Section \ref{Ex}.

\medskip
The author is grateful for helpful discussions with Tom Ilmanen and Andrew Kresch.


\section{Preliminaries from analysis}
\label{Prelim}

Consider integers $m,n\ge0$. According to one convention, a $\BC$-valued function on a closed subset $X\subset\BR^n$ lies in $C^m(X)$ if and only if it is the restriction of a $C^m$-function on an open neighborhood of~$X$. We use a different convention, following Bell \cite[page 3]{Bell1992}.

Consider an open convex subset $U\subset\BR^n$. We will say that a function $\OU\to\BC$ lies in $C^m(\OU)$ if and only if it is continuous and its restriction lies in $C^m(U)$ and all partial derivatives of order $\le m$ thereof are continuous and extend to continuous functions on~$\OU$. 

A fundamental fact from real analysis states that a function on $U$ lies in $C^m(U)$ if and only if it and all its partial derivatives of order $\le m$ exist and are continuous. A slight adaptation of the proof of this fact and the Taylor approximation yields:

\begin{Prop}\label{Taylor}
Consider any function $v\in C^m(\OU)$ and any point $\Ux\in\OU$. Then for all $\Uy\in\OU$ tending to $\Ux$ we have
$$v(\Uy)\ =\!\!\! \sum_{\underline{\nu}=(\nu_1,\ldots,\nu_n)}
\frac{1}{\nu_1!\cdots\nu_n!}\cdot
\Bigl(\frac{\partial^{\nu_1}}{\partial x_1^{\nu_1}}
\cdots\frac{\partial^{\nu_n}}{\partial x_n^{\nu_n}}f\Bigr)(\Ux) 
\cdot (\Uy-\Ux)^{\underline{\nu}} \ +\ o\bigl(|\Uy-\Ux|^m\bigr),$$
where the sum extends over all $\nu_1,\ldots,\nu_n\ge0$ with $\nu_1+\ldots+\nu_n\le m$.
\end{Prop}

From this we can directly deduce:

\begin{Prop}\label{Restriction}
For any function $v\in C^m(\OU)$ and any $C^m$-submanifold $X\subset\OU$, the restriction $v|X$ is a $C^m$-function in the sense of manifolds.
\end{Prop}

Now let $\BB$ denote the closed unit disc in~$\BC$, as before. For any function $v\in C(\BB\times\OU)$ and any point $(z,\Ux)\in\BB\times\OU$ we set
\UseTheoremCounterForNextEquation
\begin{equation}\label{KDef2}
K(v)(z,\Ux)\ := \frac{1}{2\pi i}\int_\BB \frac{v(\zeta,\Ux)}{\zeta-z}\;d\zeta\wedge\overline{d\zeta}.
\end{equation}
This integral converges, because $v$ is continuous and $\zeta\mapsto\frac{1}{\zeta-z}$ is locally $L^1$ for the measure $d\zeta\wedge\overline{d\zeta}$. Thus $K$ is a linear operator sending continuous functions on $\BB\times\OU$ to functions on $\BB\times\OU$.
Abbreviate $\bar\partial := \frac{\partial}{\partial\bar z}$, and let $\nabla$ denote the total derivative with respect to~$\Ux$.

\begin{Thm}\label{Bell2}
For any integer $m\ge1$ and any function $v\in C^{2m+1}(\BB\times\OU)$ we have $K(v)\in C^m(\BB\times\OU)$ and $\bar\partial K(v)=v$ and $\nabla K(v)=K(\nabla v)$.
\end{Thm}

\begin{Proof}
The analogous statement with $m=\infty$ and without the additional factor $\OU$ is that of Bell \cite[Thm.\;2.2]{Bell1992}. The proof given there also works with only finitely often differentiable functions and with parameters and yields the stated (perhaps suboptimal) result. Specifically set $\rho(z,\Ux) := |z|-1$. Keeping track of the order of differentiability, the same construction as in the proof of \cite[Lem.\;2.3]{Bell1992} shows that for any integers $\ell>m+1>0$ and any function $v\in C^\ell(\BB\times\OU)$ there exist functions $\Phi_m,\Psi_m\in C^{\ell-m-1}(\BB\times\OU)$ such that $\Phi_m|(\partial\BB\times\OU)=0$ and $v=\bar\partial\Phi_m+\Psi_m\rho^{m+1}$. Applying this with $\ell=2m+1$, the rest of the proof goes through likewise.
\end{Proof}

\section{Proof of the strong Nullstellensatz}
\label{NSProof}

In this section we fix an ideal $I\subset\BC[\UX]$ and choose a system of generators $(f_1,\ldots,f_m)$. Consider the real analytic function 
\UseTheoremCounterForNextEquation
\begin{equation}\label{FDef}
F\colon \BC^n\to\BR^{\ge0},\ 
\textstyle \Uz\mapsto\sum_{j=1}^m |f_j(\Uz)|^2.
\end{equation}
Its zero locus on $\BB^n$ is precisely the subset $M(I)$ from (\ref{MIDef}). Fix a polynomial $g\in I(M(I))$. Then by construction, the zero locus of $F$ on $\BB^n$ is contained in the zero locus of~$g$. By the {\L}ojasiewicz inequality (see Bierstone-Milman \cite[Thm.\;6.4]{BierstoneMilman1988}) we can therefore choose an integer $r>0$ and a real number $c>0$ such that
\UseTheoremCounterForNextEquation
\begin{equation}\label{Loja}
\forall\Uz\in\BB^n\colon \ |g(\Uz)|^r \le c\cdot F(\Uz).
\end{equation}
In the rest of this section we identify polynomials in $\BC[\UX]$ with the functions $\BB^n\to\BC$ that they represent. For any integers $s,t>0$ consider the function
\UseTheoremCounterForNextEquation
\begin{equation}\label{phiDef}
\phi_{s,t}\colon \BB^n\to\BC,\ \Uz\mapsto 
\left\{\vcenter{\hrule width0pt height40pt}\right.\!\!
\begin{array}{cl}
\displaystyle\frac{g(\Uz)^t}{F(\Uz)^s} & \hbox{if $\Uz\not\in M(I)$,} \\[13pt]
0 & \hbox{if $\Uz\in M(I)$.}
\end{array}
\end{equation}

\begin{Lem}\label{phiLem}
For any integers $\ell\ge0$ and $s>0$ and $t>r(s+\ell)$ we have $\phi_{s,t}\in C^\ell(\BB^n)$.
\end{Lem}

\begin{Proof}
By construction $\phi_{s,t}$ is $C^\infty$ outside $M(I)$. If $t>rs$, then on $\BB\setminus M(I)$ we have
\UseTheoremCounterForNextEquation
\begin{equation}\label{phiLema}
\bigl|\phi_{s,t}\bigr|
\ =\ \frac{|g|^t}{F^s}
\ =\ \biggl(\frac{|g|^r}{F}\biggr)^s\cdot|g|^{t-rs}
\ \stackrel{(\ref{Loja})}{\le} \ c^s\cdot|g|^{t-rs}.
\end{equation}
Here $|g|^{t-rs}$ is continuous and vanishes on $M(I)$, hence the same holds for~$\phi_{s,t}$. This proves the lemma in the case $\ell=0$. 

Next assume that $t>r(s+1)$. Then $t-rs\ge2$, and so (\ref{phiLema}) implies that $\phi_{s,t}$ is real differentiable everywhere and its total derivative $\nabla\phi_{s,t}$ vanishes along $M(I)$. On $\BB^n\setminus M(I)$ we have
$$\nabla\phi_{s,t}
\ =\ t\cdot\frac{g^{t-1}}{F^s}\cdot\nabla g
          - s\cdot\frac{g^t}{F^{s+1}}\cdot\nabla F$$
and hence
\UseTheoremCounterForNextEquation
\begin{equation}\label{phiLemb}
\nabla\phi_{s,t}
\ =\ t\cdot\phi_{s,t-1}\cdot\nabla g
          - s\cdot\phi_{s+1,t}\cdot\nabla F.
\end{equation}
Applying (\ref{phiLema}) to $(s,t-1)$ and $(s+1,t)$ in place of $(s,t)$ shows that $\phi_{s,t-1}$ and $\phi_{s+1,t}$ are continuous and vanish on $M(I$). Thus the equation (\ref{phiLemb}) holds on all of~$\BB^n$, and $\phi_{s,t}$ is $C^1$. This proves the lemma in the case $\ell=1$. 

Equation (\ref{phiLemb}) now also implies the general case $\ell\ge1$ by induction on~$\ell$.
\end{Proof}

\medskip
The next lemma and its proof are adapted from Krantz-Li \cite[\S2]{KrantzLi1996}. For all $1\le i\le n$ we abbreviate $\bar\partial_i := \frac{\partial}{\partial\bar z_i}$ 

\begin{Lem}\label{gLem}
Consider any integers $\ell\ge3$ and $t>r(2^n\ell-1)$. Then for every integer ${0\le k\le n}$ there exist functions $g_{j,k} \in C^{2^{n-k}\ell-2}(\BB^n)$ for all $1\le j\le m$ satisfying $\bar\partial_ig_{j,k}=0$ for all $1\le i\le k$, such that
$$g^{2^kt}\ =\ \sum_{j=1}^m g_{j,k}f_j.$$
\end{Lem}

\begin{Proof}
We prove this by induction on~$k$, beginning with $k=0$. By the construction (\ref{phiDef}) with $s=1$ and by (\ref{FDef}) we have
$$g^t\ =\ \phi_{1,t}F\ =\ 
\phi_{1,t}\sum_{j=1}^m \overline{f_j}f_j
\ =\ \sum_{j=1}^m \phi_{1,t}\overline{f_j}\cdot f_j$$
on $\BB^n$. Here $\phi_{1,t}\in C^{2^n\ell-2}(\BB^n)$ by Lemma \ref{phiLem}, and $\overline{f_j}$ is already $C^\infty$ everywhere. Thus the functions $g_{j,0} := \phi_{1,t}\overline{f_j}$ possess the desired properties for $k=0$. 

Now assume that the desired functions $g_{j,k}$ are already given for some fixed $0\le k<n$. We must construct the next batch of functions $g_{j,k+1}$. For any function $h\in C(\BB^n)$ we set
$$K_{k+1}(h)(z_1,\ldots,z_n)
\ := \frac{1}{2\pi i}\int_\BB \frac{f(z_1,\ldots,z_k,\zeta,z_{k+2},\ldots,z_n)}{\zeta-z_{k+1}}\;d\zeta\wedge\overline{d\zeta}.$$
Up to reordering the variables the operator $K_{k+1}$ is just the operator $K$ from (\ref{KDef2}). 
For any indices $1\le j,j'\le m$ we have $g_{j',k}$, $g_{j,k}\in C^{2^{n-k}\ell-2}(\BB^n)$ by the induction hypothesis, and hence 
$g_{j',k}\,\bar\partial_{k+1}g_{j,k}\in C^{2^{n-k}\ell-3}(\BB^n)$. By Theorem \ref{Bell2} we therefore have
$$u_{j,j'}\ :=\ K_{k+1}\bigl(g_{j',k}\,\bar\partial_{k+1}g_{j,k}\bigr)
\ \in\ C^{2^{n-k-1}\ell-2}(\BB^n)$$
and hence
$$g_{j,k+1}\ := g^{2^kt}g_{j,k} - 
\sum_{j'=1}^m (u_{j,j'}-u_{j',j})\, f_{j'}
\ \in\ C^{2^{n-k-1}\ell-2}(\BB^n).$$
The definition of $g_{j,k+1}$ implies that
\begin{eqnarray*}
\sum_{j=1}^m g_{j,k+1}f_j
&\!\!=\!\!& \sum_{j=1}^m g^{2^kt}\,g_{j,k}f_j 
- \sum_{j,j'=1}^m (u_{j,j'}-u_{j',j}) f_{j'}f_j \\
&\!\!=\!\!& g^{2^kt}\sum_{j=1}^m g_{j,k}f_j 
- \sum_{j,j'=1}^m u_{j,j'} f_{j'}f_j 
+ \sum_{j,j'=1}^m u_{j',j} f_{j'}f_j \\
&\!\!=\!\!& g^{2^kt}\, g^{2^kt}\ =\ g^{2^{k+1}t},
\end{eqnarray*}
using the induction hypothesis and the symmetry between the last two sums.

Next, for all $1\le i\le k$ we have $\bar\partial_i g_{j',k}=\bar\partial_i g_{j,k}=0$ by the induction hypothesis, and since $g_{j,k}\in C^2$, therefore also $\bar\partial_i (\bar\partial_{k+1} g_{j,k}) = \bar\partial_{k+1} (\bar\partial_i g_{j,k}) = 0$. With the last equation in Theorem \ref{Bell2} we deduce that
$$\bar\partial_i u_{j,j'}
\ =\ \bar\partial_i K_{k+1}\bigl(g_{j',k}\,\bar\partial_{k+1}g_{j,k}\bigr)
\ =\ K_{k+1}\bigl(\bar\partial_i\bigl(g_{j',k}\,\bar\partial_{k+1}g_{j,k}\bigr)\bigr)
\ =\ K_{k+1}(0)\ =\ 0.$$
Since $f_{j'}$ and $g$ are holomorphic, plugging this into the formula defining $g_{j,k+1}$ and using the induction hypothesis shows that
$$\bar\partial_i g_{j,k+1}
\ =\ g^{2^kt} \; \bar\partial_i g_{j,k} \ =\ 0.$$
Moreover, by the definition of $u_{j,j'}$ and Theorem \ref{Bell2} we have
$$\bar\partial_{k+1} u_{j,j'}\ =\ \bar\partial_{k+1} K_{k+1}\bigl(g_{j',k}\,\bar\partial_{k+1}g_{j,k}\bigr)
\ =\ g_{j',k}\,\bar\partial_{k+1}g_{j,k}.$$
Using this and the holomorphy of $f_{j'}$ and $g$ we calculate
\begin{eqnarray*}
\bar\partial_{k+1} g_{j,k+1}
&\!\!=\!\!& g^{2^kt} \; \bar\partial_{k+1} g_{j,k} - 
\sum_{j'=1}^m \bigl(\bar\partial_{k+1} u_{j,j'}-\bar\partial_{k+1} u_{j',j}\bigr)\, f_{j'} \\
&\!\!=\!\!& g^{2^kt} \; \bar\partial_{k+1} g_{j,k} - 
\sum_{j'=1}^m \bigl(g_{j',k}\,\bar\partial_{k+1}g_{j,k}\!- g_{j,k}\,\bar\partial_{k+1}g_{j',k}\bigr)\, f_{j'} \\
&\!\!=\!\!& 
\biggl( g^{2^kt} - \sum_{j'=1}^m g_{j',k}f_{j'}\biggr)\,
\bar\partial_{k+1} g_{j,k}
+ \; g_{j,k}\, \bar\partial_{k+1}\biggl(\;\sum_{j'=1}^m g_{j',k}\cdot f_{j'}\biggr).
\end{eqnarray*}
With the induction hypothesis and the holomorphy of $g$ we conclude that
$$\bar\partial_{k+1} g_{j,k+1}
\ =\ (g^{2^kt}-g^{2^kt})\,\bar\partial_{k+1} g_{j,k}
+ \; g_{j,k}\,\bar\partial_{k+1}(g^{2^kt})\ =\ 0.$$
Thus the functions $g_{j,k+1}$ satisfy all requirements, and the lemma is proved.
\end{Proof}

\medskip
The case $k=n$ of Lemma \ref{gLem} directly yields:

\begin{Thm}\label{NS1a}
For any ideal $I=(f_1,\ldots,f_m)\subset\BC[\UX]$, any polynomial $g\in I(M(I))$, and any integer $\ell\ge3$ there exist an integer $N\ge1$ and functions $g_j\in C^{\ell-2}(\BB^n)$ which are holomorphic on $(\BB^\circ)^n$, such that
$$g^N\ =\ \sum_{j=1}^m g_j f_j.$$
\end{Thm}

Now let $A$ denote the ring of all power series in $\BC[[\UX]]$ satisfying
\UseTheoremCounterForNextEquation
\begin{equation}\label{NormDefa}
\textstyle \bigl\Vert\sum_\Uk a_\Uk \UX^\Uk\bigr\Vert \ :=\ \sum_\Uk |a_\Uk| \ <\ \infty
\end{equation}
Equivalently this is the ring of all power series that converge absolutely on~$\BB^n$, or again the Banach algebra completion of $\BC[\UX]$ with respect to $\Vert\ \Vert$.

\begin{Lem}\label{AbsConvLem}
For any integer $k>\frac{n}{2}$, any function $h\in C^k(\BB^n)$ which is holomorphic on $(\BB^\circ)^n$ is represented by a power series in~$A$.
\end{Lem}

\begin{Proof}
Consider the Fourier series of $h|(\partial\BB)^n$, written as a Laurent series $\sum_{\Uk\in\BZ^n} a_\Uk \UX^\Uk$ with 
$$a_\Uk\ :=\ \biggl(\frac{1}{2\pi i}\biggr)^n \int_{(\partial\BB)^n} \frac{h(\Uz)}{\Uz^\Uk}\, \frac{dz_1}{z_1}\cdots\frac{dz_n}{z_n}.$$
Since $h$ is continuous, this coefficient is the limit for $r\nearrow1$ of
$$a_\Uk(r)\ :=\ \biggl(\frac{1}{2\pi i}\biggr)^n \int_{(\partial\BB)^n} \frac{h(r\Uz)}{(r\Uz)^\Uk}\, \frac{d(rz_1)}{rz_1}\cdots\frac{d(rz_n)}{rz_n}.$$
As $h|(\BB^\circ)^n$ is holomorphic, by the Cauchy integral formula $a_\Uk(r)$ is zero unless $\Uk\in(\BZ^{\ge0})^n$, in which case it is the coefficient of $\UX^\Uk$ in the power series representing $h|(\BB^\circ)^n$. In particular $a_\Uk(r)$ is independent of $r$ and hence equal to~$a_\Uk$. Thus $\sum_\Uk a_\Uk \UX^\Uk$ is really the power series representing $h|(\BB^\circ)^n$. 

On the other hand the restriction $h|(\partial\BB)^n$ is $C^k$ by the assumption and Proposition \ref{Restriction}. Since $k>\frac{n}{2}$, its Fourier series is therefore absolutely convergent (see Grafakos \cite[Thm.\,3.3.16]{Grafakos2014}). Thus the power series $\sum_\Uk a_\Uk \UX^\Uk$ is absolutely convergent on~$\BB^n$. Finally, since it represents $h$ on the interior $(\BB^\circ)^n$, by continuity it represents $h$ on all of $\BB^n$.
\end{Proof}

\begin{Thm}\label{NS1b}
For any ideal $I=(f_1,\ldots,f_m)\subset\BC[\UX]$ and any polynomial $g\in I(M(I))$, there exist an integer $N\ge1$ and power series $g_j\in A$ such that
$$g^N\ =\ \sum_{j=1}^m g_j f_j.$$
\end{Thm}

\begin{Proof}
Use Theorem \ref{NS1a} with any $\ell\ge\frac{n}{2}+3$ and apply Lemma \ref{AbsConvLem}.
\end{Proof}

\begin{Thm}\label{NS1rep}
{\bf(=\ \ref{NS1})}\ 
For any ideal $I\subset\BC[\UX]$ we have $I(M(I)) = \Rad(\OI)$.
\end{Thm}

\begin{Proof}
As before write $I=(f_1,\ldots,f_m)$. For any $g\in I(M(I))$ choose $N$ and $g_j$ as in Theorem \ref{NS1b}. Write each $g_j$ as the limit in $A$ of a sequence of polynomials $g_{j,k}\in\BC[\UX]$ for $k\to\infty$. Then $g^N=\sum_{j=1}^m g_j f_j$ is the limit in $A$ of the sequence of polynomials $\sum_{j=1}^m g_{j,k} f_j \in I$ for $k\to\infty$. Since $g^N$ is already a polynomial, this limit process already takes place in the normed algebra $\BC[\UX]$; hence $g^N\in\OI$. This proves that $I(M(I)) \subset \Rad(\OI)$.

Conversely consider any point $\Uz\in M(I)$. Then by (\ref{MIDef}) we have $f(\Uz)=0$ for all $f\in I$. Since evaluation at $\Uz$ defines a continuous map $\BC[\UX]\to\BC$ for the norm $\Vert\ \Vert$, it follows that $f(\Uz)=0$ for all $f\in\OI$ as well. By the definition of the radical the same then also follows for all $f\in\Rad(\OI)$. Varying $\Uz\in M(I)$ and using (\ref{IZDef}) thus shows that $\Rad(\OI)\subset I(M(I))$. 
\end{Proof}


\section{Ideals with finite support}
\label{FinSupp}

We begin with the case of one variable~$X$.

\begin{Lem}\label{1var1}
For any $z\in\BC$ with $|z|>1$ we have $1\in\overline{(X-z)}$.
\end{Lem}

\begin{Proof}
For all $m\ge1$ we have $1-\bigl(\frac{X}{z}\bigr)^m\in(X-z)$. Since $|z|>1$, we have $\bigl\Vert\bigl(\frac{X}{z}\bigr)^m\bigr\Vert = \frac{1}{z^m}\to0$ for $m\to\infty$. In the limit we deduce that $1\in\overline{(X-z)}$.
\end{Proof}


\begin{Lem}\label{1var3}
For any $z\in\BC$ with $|z|=1$ and any $k\ge1$ we have $(X-z)\in\overline{((X-z)^k)}$.
\end{Lem}

\begin{Proof}
For all $m\ge1$ the binomial theorem shows that
$$X^m\ =\ (z+(X-z))^m\ \equiv\ z^m+mz^{m-1}(X-z)
\ \ \hbox{modulo}\ \ ((X-z)^2).$$
Equivalently
$$(X-z)\ \equiv\ \frac{1}{mz^{m-1}}(X^m-z^m)\ \ \hbox{modulo}\ \ ((X-z)^2).$$
Since $|z|=1$, we have $\bigl\Vert\frac{1}{mz^{m-1}}(X^m-z^m)\bigr\Vert = \frac{1}{m}(1+1)\to0$ for $m\to\infty$. In the limit we deduce that $(X-z)\equiv0$ modulo $\overline{((X-z)^2)}$. This shows the case $k=2$ of the lemma. 

The general case follows by induction on~$k$, the case $k=1$ being trivial. If $k>2$, the induction hypothesis implies that $(X-z)^2 \in (X-z)\cdot \overline{((X-z)^{k-1})} \subset \overline{((X-z)^k)}$. Using the case $k=2$ it follows that $(X-z) \in \overline{((X-z)^2)} \subset \overline{((X-z)^k)}$, as desired.
\end{Proof}

\begin{Lem}\label{1var4}
For any $z\in\BC$ with $|z|=1$ the ideal $(X-z)$ is closed.
\end{Lem}

\begin{Proof}
Since $|z|=1$, the evaluation map $\ell\colon\BC[X]\to\BC$, $f\mapsto f(z)$ satisfies $|\ell(f)|\le\Vert f\Vert$. It is therefore continuous with respect to the metric induced by~$\Vert\ \Vert$, and so its kernel $(X-z)$ is closed.
\end{Proof}

\begin{Lem}\label{1var5}
For any $z\in\BC$ with $|z|<1$ and any $k\ge1$ the ideal $((X-z)^k)$ is closed.
\end{Lem}

\begin{Proof}
For any fixed $\nu\ge0$ consider the linear map $\ell_\nu\colon\BC[X]\to\BC$, $f\mapsto f^{(\nu)}(z)$. For each $j\ge0$ we have $\ell_\nu(X^j) = j(j-1)\cdots(j-\nu+1)z^{j-\nu}$. Since $|z|<1$, this value tends to $0$ for $j\to\infty$. Thus there exists a real number $c_\nu>0$ such that $|\ell_\nu(X^j)|\le c$ for all $j\ge0$. It then follows that $|\ell_\nu(f)|\le c_\nu\Vert f\Vert$ for all $f\in\BC[X]$. Therefore $\ell_\nu$ is continuous with respect to~$\Vert\ \Vert$, and so its kernel is a closed subspace.

Varying $\nu$ it now follows that $((X-z)^k) = \bigcap_{\nu=0}^{k-1}\Ker(\ell_\nu)$ is closed.
\end{Proof}

\medskip
Now we return to an arbitrary number of variables. 

\begin{Lem}\label{FiniteCodimLem}
For any ideals $J\subset I\subset\BC[\UX]$ with $\dim_\BC(I/J)<\infty$, if $J$ is closed, so is~$I$.
\end{Lem}

\begin{Proof}
Consider the seminorm induced by $\Vert\ \Vert$ on the factor space $\BC[\UX]/J$. Since $J$ is closed, this seminorm is a norm. As any finite dimensional subspace of a normed $\BC$-vector space is closed, it follows that $I/J\subset\BC[\UX]/J$ is closed for the induced norm. Its inverse image $I\subset\BC[\UX]$ under the projection map is therefore closed for the norm $\Vert\ \Vert$.
\end{Proof}

\medskip
For any point $\Uz=(z_1,\ldots,z_n)\in\BC^n$ consider the maximal ideal 
$$\Fm_\Uz\ :=\ (X_1-z_1,\ldots,{X_n-z_n})\ \subset\ \BC[\UX].$$

\begin{Thm}\label{NS2rep}
{\bf(=\ \ref{NS2})}\ 
For any point $\Uz\in\BC^n$ and any $\Fm_\Uz$-primary ideal $I\subset\BC[\UX]$ we have
$$\OI\ =\ \
\biggl\{\!\begin{array}{ll}
I + (X_j-z_j)_{1\le j\le n,\, |z_j|=1}
 & \hbox{if $\Uz\in\BB^n$,} \\[3pt]
(1) & \hbox{if $\Uz\not\in\BB^n$.}
\end{array}$$
\end{Thm}

\begin{Proof}
We deduce this from the one variable case using the isometric embeddings $\BC[X_j]\into\BC[\UX]$ for all $1\le j\le n$. By assumption we have $\Fm_\Uz^r\subset I$ for some $r\ge1$. 

If $|z_j|>1$ for some~$j$, we have $1\in\overline{(X_j-z_j)}$ by Lemma \ref{1var1}, and hence $1\in\bigl(\,\overline{(X_j-z_j)}\,\bigr)^r \allowbreak \subset \overline{((X_j-z_j)^r)} \subset \OI$. Thus $\OI=(1)$, as desired. 

So suppose that $|z_j|\le1$ for all $1\le j\le n$. By symmetry we can assume that $|z_j|<1$ if $j\le m$, and $|z_j|=1$ if $j>m$. 
Since $(X_j-z_j)^r\in \Fm_\Uz^r\subset I$, for each $m<j\le n$ we have $(X_j-z_j)\in \overline{((X_j-z_j)^r)} \subset \OI$ by Lemma \ref{1var3}. Replacing $I$ by $I + (X_j-z_j)_{m<j\le n}$ therefore does not change~$\OI$. Then we have $J\subset I \subset \BC[\UX]$ with the ideal
$$J\ :=\ \bigl((X_1-z_1)^r,\ldots,(X_m-z_m)^r,(X_{m+1}-z_{m+1}),\ldots,(X_n-z_n)\bigr).$$
It remains to show that any ideal $I$ with this property is closed. For this observe that for any indices $\nu_1,\ldots,\nu_m\ge0$, the linear map 
$$\BC[\UX]\to\BC,\ f\mapsto \bigl(\tfrac{\partial^{\nu_1}}{\partial X_1^{\nu_1}} \cdots\tfrac{\partial^{\nu_m}}{\partial X_m^{\nu_m}}f\bigr)(\Uz)$$
is continuous with respect to $\Vert\ \Vert$ by the same arguments as in the proofs of Lemmas \ref{1var4} and \ref{1var5}. Thus its kernel is closed. Since $J$ is the intersection of these kernels for all possible indices $\nu_1,\ldots,\nu_m\in\{0,\ldots,r-1\}$, it follows that $J$ is closed. With Lemma \ref{FiniteCodimLem} we deduce that $I$ is closed, as desired.
\end{Proof}

\medskip
Note that Theorem \ref{NS2} yields an explicit description of $\overline{\Fm_\Uz^r}$ for all $\Uz\in\BB^n$ and $r\ge1$. 

\begin{Def}\label{ITildeDef}
To any ideal $I\subset\BC[\UX]$ we associate the ideal
$$\TI\ :=\ \bigcap_{\Uz\in\BB^n} \bigcap_{r\ge1} \bigl(I+\overline{\Fm_\Uz^r}\,\bigr).$$
\end{Def}

\begin{Prop}\label{ITildeProp}
\begin{enumerate}
\item[(a)] $\TI$ is closed.
\item[(b)] $I\subset\TI\subset I(M(I))$.
\item[(c)] $\OI\subset\TI=\widetilde{\OI}=\widetilde{\TI}$.
\item[(d)] For any ideals $I\subset J$ we have $\TI\subset\TJ$.
\item[(e)] For any ideals $I,J$ we have $\TI\cdot\TJ\subset\widetilde{(IJ)}$.
\end{enumerate}
\end{Prop}

\begin{Proof}
The ideals $I+\overline{\Fm_\Uz^r}$ are all closed by Lemma \ref{FiniteCodimLem}. Thus their intersection is closed, proving (a). 
Next, the inclusion $I\subset\smash{\TI}$ is obvious. Also, for all $\Uz\in M(I)$ we have $I\subset\Fm_\Uz$ and $\overline{\Fm_\Uz}=\Fm_\Uz$ and hence ${\smash\TI}\subset I+\overline{\Fm_\Uz}=\Fm_\Uz$. Therefore ${\smash\TI}\subset\bigcap_{\Uz\in M(I)}\Fm_\Uz = I(M(I))$, proving (b).
Assertion (d) is a direct consequence of the definition. 

Returning to (c), for all $\Uz\in\BB^n$ and $r\ge1$ the definition implies that $\TI+\overline{\Fm_\Uz^r} \subset \bigl(I+\overline{\Fm_\Uz^r}\,\bigr)+\overline{\Fm_\Uz^r} = I+\overline{\Fm_\Uz^r}$. Varying $\Uz$ and $r$ it follows that $\widetilde{\TI}\subset\TI$. On the other hand, by (a) we have $I\subset\OI\subset\TI$, which by (d) implies that
$\TI\subset\widetilde{\OI}\subset\widetilde{\TI}$. Together this implies the equalities in (d).

Finally in (e), for all $\Uz\in\BB^n$ and $r\ge1$ we have $\TI\cdot\TJ \subset \bigl(I+\overline{\Fm_\Uz^r}\,\bigr) \cdot\bigl(J+\overline{\Fm_\Uz^r}\,\bigr) \subset \bigl(IJ+\overline{\Fm_\Uz^r}\,\bigr)$. By varying $\Uz$ and $r$ this implies that $\TI\cdot\TJ\subset\widetilde{(IJ)}$, as desired.
\end{Proof}

\begin{Thm}\label{NS4a}
For any ideal $I\subset\BC[\UX]$ with $M(I)$ finite we have $\OI=\TI$.
\end{Thm}

\begin{Proof}
By Proposition \ref{ITildeProp} (c) we already have $\OI\subset\TI$. Conversely, by the finiteness of $M(I)$ and Theorem \ref{NS1} we have $\Rad(\OI)=\prod_{\Uz\in M(I)} \Fm_\Uz$.
Thus there exists $r\ge1$ with $\prod_{\Uz\in M(I)} \Fm_\Uz^r \subset \OI$. By continuity of addition and multiplication in $\BC[\UX]$ this implies that $\prod_{\Uz\in M(I)} \overline{\Fm_\Uz^r} \subset \OI$. 
As the ideals $I+\overline{\Fm_\Uz^r}$ are primary to mutually distinct maximal ideals, it follows that
$$\TI\ \subset\ \bigcap_{\Uz\in M(I)} \bigl(I+\overline{\Fm_\Uz^r}\,\bigr)
\ \stackrel{!}{=} \prod_{\Uz\in M(I)} \bigl(I+\overline{\Fm_\Uz^r}\,\bigr)
\ \subset\ \OI.$$
Together this implies that $\OI=\TI$.
\end{Proof}


\section{Plane curves}
\label{PlaneCurves}

In the rest of this article we study the case $n=2$. 

\begin{Prop}\label{PlaneClass}
For any irreducible polynomial $f\in\BC[X_1,X_2]$ we have precisely one of the following cases:
\begin{enumerate}
\item[(a)] $M((f))\cap(\BB^\circ)^2\not=\emptyset$.
\item[(b)] $M((f))=\{z_1\}\times\BB$ and $f=u(X_1-z_1)$ for some $z_1\in\partial\BB$ and $u\in\BC^\times$.
\item[(c)] $M((f))=\BB\times\{z_2\}$ and $f=u(X_2-z_2)$ for some $z_2\in\partial\BB$ and $u\in\BC^\times$.
\item[(d)] $M((f))$ is infinite and contained in $(\partial\BB)^2$.
\item[(e)] $M((f))$ is finite and contained in $(\partial\BB)^2$, possibly empty.
\end{enumerate}
Moreover, in the cases (a) through (d) we have $I(M((f)))=(f)$.
\end{Prop}

\begin{Proof}
If $M((f))$ contains a point $\Uz\in(\BB^\circ)^2$, it contains a whole neighborhood of $\Uz$ in the irreducible curve $V((f))$. Then $M((f))$ is Zariski dense in $V((f))$; hence $I(M((f)))=(f)$, and we have the case (a).

Next suppose that $M((f))$ contains a point $\Uz=(z_1,z_2)\in\partial\BB\times\BB^\circ$. If the projection to the first coordinate $V((f))\to\BC$ is not constant, it is an open map, and so any neighborhood of $\Uz$ contains a point $(z_1',z_2')$ with $z_1'\in\BB^\circ$. Choosing the neighborhood small enough guarantees that $z_2'\in\BB^\circ$ as well, and we are back in the case (a). Otherwise the projection map is constant and the curve must be given by the equation $X_1=z_1$. Then $M((f))=\{z_1\}\times\BB$ and of course $I(M((f)))=(f)$, and we have the case (b).

By symmetry, if $M((f))$ contains a point from $\BB^\circ\times\partial\BB$, we have the case (a) or (c).

If none of these cases applies, we must have $M((f))\subset(\partial\BB)^2$. If $M((f))$ is then infinite, it is again Zariski dense in the irreducible curve $V((f))$; hence $I(M((f)))=(f)$, and we have the case (d). Otherwise we are left with the case (e).
\end{Proof}

\begin{Prop}\label{PlaneA}
In the case (a) of Proposition \ref{PlaneClass}, for all $k\ge1$ we have 
$$\overline{(f^k)} \ =\ \widetilde{(f^k)} \ =\ (f^k).$$
\end{Prop}

\begin{Proof}
Pick any point $\Uz\in M((f))\cap(\BB^\circ)^2$. Then for each $r\ge1$ we have $\overline{\Fm_\Uz^r}=\Fm_\Uz^r$ by Theorem \ref{NS2}. As the completion of $\BC[\UX]$ for the $\Fm_\Uz$-adic topology is the power series ring $\BC[[\UX-\Uz]] := \BC[[X_1-z_1,X_2-z_2]]$, we deduce that
$$\widetilde{(f^k)}\ \subset\ I \ :=\ \bigcap_{r\ge1} \bigl((f^k)+\Fm_\Uz^r\,\bigr)\ =\ \BC[\UX] \cap f^k\cdot\BC[[\UX-\Uz]].$$
We claim that $I=(f^k)$. This, together with the inclusions $(f^k) \subset \overline{(f^k)} \subset \widetilde{(f^k)}$ from Proposition \ref{ITildeProp} yields the desired equalities.

To prove the claim observe first that, since $f$ lies in $\Fm_\Uz$, it is not a unit in $\BC[[\UX-\Uz]]$. Thus $f\cdot\BC[[\UX-\Uz]]$ is an ideal of height $1$ of $\BC[[\UX-\Uz]]$, and so $\dim_\BC\bigl(\BC[[\UX-\Uz]]/f\cdot\BC[[\UX-\Uz]]\bigr)=\infty$. As the image of $\BC[\UX]$ in $\BC[[\UX-\Uz]]/f\cdot\BC[[\UX-\Uz]]$ is dense for the $(X_1-z_1,X_2-z_2)$-adic topology, it, too, is infinite dimensional. Thus $\BC[\UX] \cap f\cdot\BC[[\UX-\Uz]]$ is an ideal of infinite codimension of $\BC[\UX]$ containing the irreducible polynomial~$f$. This ideal is therefore equal to $(f)$. By induction on $k$ it follows that $\BC[\UX] \cap f^k\cdot\BC[[\UX-\Uz]] = (f^k)$ for all $k\ge1$, as claimed.
\end{Proof}

\begin{Prop}\label{PlaneBC}
In the cases (b) and (c) of Proposition \ref{PlaneClass}, for all $k\ge1$ we have 
$$\overline{(f^k)} \ =\ \widetilde{(f^k)} \ =\ (f).$$
\end{Prop}

\begin{Proof}
By symmetry it suffices to consider the case \ref{PlaneClass} (b). Without loss of generality we then have $f=(X_1-z_1)$ for some $z_1\in\partial\BB$. By Lemma \ref{1var3} and the isometric embedding $\BC[X_1]\into\BC[\UX]$ we already know that $(f)\subset\overline{(f^k)}$. On the other hand Proposition \ref{ITildeProp} (c) and (b) implies that $\overline{(f^k)} \subset \widetilde{(f^k)} \subset (f)$. Together this yields the desired equalities.
\end{Proof}

\medskip
In the case (e) of Proposition \ref{PlaneClass} the ideal $\overline{(f^k)}$ is described by Theorem \ref{NS4a} for all $k\ge1$.
The most interesting and difficult case (d) of Proposition \ref{PlaneClass} is treated in the following section.


\section{Plane curves touching the bidisc}
\label{PlaneD}

This section is devoted to the case (d) of Proposition \ref{PlaneClass}. First we give an explicit description of the curves with this property.

Consider an irreducible polynomial $f\in\BC[X_1,X_2]$ for which $M((f))$ is infinite and contained in $(\partial\BB)^2$. We view $\BC$ as a subset of the Riemann sphere $\hat\BC := \BC\cup\{\infty\} \cong\BP^1(\BC)$, and let $C\subset\hat\BC^2$ denote the closure of the curve $V(f)$. Let $\pi\colon\tilde C\to C$ denote the normalization of~$C$, so that $\tilde C$ is an irreducible smooth projective algebraic curve over~$\BC$. Let $\pi_j\colon \tilde C\to\hat\BC$ denote the composite of $\pi$ with the projection to the $j$-th factor.
Recall that the M\"obius transformation $\mu(z) := i\frac{1-z}{1+z}$ is an automorphism of $\hat\BC$ with 
$$\begin{array}{rcl}
\mu(\partial\BB)&\!=\ \hat\BR \ :=\!& \BR\cup\{\infty\} \ \cong\ \BP^1(\BR)
\ \ \hbox{and} \\[3pt]
\mu(\BB^\circ)  &\!=\ \BH \ :=\!& \{z\in\BC\mid \Im(z)>0\}.
\end{array}$$

\begin{Prop}\label{PlaneDClass}
\begin{enumerate}
\item[(a)] The curve $(\mu\times\mu)(C)\subset\hat\BC^2$ is defined over~$\BR$.
\item[(b)] We have $C\cap(\partial\BB\times\hat\BC) = C\cap(\hat\BC\times\partial\BB) = C\cap(\partial\BB)^2$.
\item[(c)] The bidegree $(d_1,d_2)$ of $C\subset\hat\BC^2$ satisfies $d_1,d_2\ge1$.
\item[(d)] For each $j$ the map $\pi_j\colon \tilde C\to\hat\BC$ is unramified over $\partial\BB$.
\item[(e)] For each $z\in\partial\BB$ we have $|\pi_1^{-1}(z)|=d_2$
and $|\pi_2^{-1}(z)|=d_1$.
\end{enumerate}
\end{Prop}

\begin{Proof}
By construction $C' := (\mu\times\mu)(C)$ is an irreducible curve in $\hat\BC^2$ with infinitely many points in $\hat\BR^2$; hence it is defined over~$\BR$, proving (a). 
Also, by assumption we have $C\cap(\partial\BB\times\BB^\circ)=\emptyset$, which is equivalent to $C'\cap(\hat\BR\times\BH)=\emptyset$. 
Since $C'$ is defined over~$\BR$, by complex conjugation it follows that $C'\cap(\hat\BR\times(-\BH))=\emptyset$. Thus $C'\cap(\hat\BR\times\hat\BC)\subset\hat\BR^2$ and hence $C\cap(\partial\BB\times\hat\BC)\subset(\partial\BB)^2$. By symmetry we obtain (b).

Next, if one coefficient of the bidegree of $C$ were $0$, the curve would have the form $\{z_1\}\times\hat\BC$ or $\hat\BC\times\{z_2\}$ and we would have the case (b) or (c) of Proposition \ref{PlaneClass}. This shows (c).

Now consider a point $\tilde c\in\tilde C$ with image $(z_1,z_2):=\pi(\tilde c)\in\hat\BC^2$. By (b) we have $z_1\in\partial\BB$ if and only if $z_2\in\partial\BB$. Assume this to be the case. After substituting each $X_j$ by $z_jX_j$, we may without loss of generality assume that each $z_j=1$. Then $\mu(\pi_j(\tilde c))=0$.
By (c) the map $\pi_j$ is non-constant of degree $d_{3-j}$. Let $e_j\ge1$ denote its ramification degree at~$\tilde c$. Choose a local chart of $\tilde C$ at $\tilde c$ with parameter~$z$, such that $\mu(\pi_1(z)) = z^{e_1}$. Then locally at~$\tilde c$, the inverse image $\pi_1^{-1}(\partial\BB)=(\mu\circ\pi_1)^{-1}(\hat\BR)$ consists of the $2e_1$ rays $\zeta\BR^{\ge0}$ for all $\zeta\in\BC$ with $\zeta^{2e_1}=1$. 
Likewise, locally near $\tilde c$ we have $\mu(\pi_2(z)) = z^{e_2}u(z)$, where $u$ is analytic at $0$ with $u(0)\not=0$. The inverse image $\pi_2^{-1}(\partial\BB)=(\mu\circ\pi_2)^{-1}(\hat\BR)$ thus consists of $2e_2$ 
smooth curve segments 
emanating from~$\tilde c$. But by (b) the two inverse images coincide. Thus $e_1$ and $e_2$ are equal, say to $e\ge1$.

By the equality of the inverse images, we can now say that there exists $\epsilon>0$ such that for all $\zeta\in\BC$ with $\zeta^{2e}=1$ and all real numbers $t\in\left]0,\epsilon\right[$ we have $\mu(\pi_2(\zeta t)) = (\zeta t)^eu(\zeta t) \in\BR$. Since $\zeta^e=\pm1$, this is equivalent to $u(\zeta t) \in\BR$. Consider the power series expansion $u(z)=\sum_{k\ge0}u_kz^k$. Then the case $\zeta=1$ and the identity theorem for power series implies that all $u_k\in\BR$. For arbitrary $\zeta$ the condition is thus equivalent to $\sum_{k\ge0}u_k(\zeta^k-\zeta^{-k})t^k = u(\zeta t)-\overline{u(\zeta t)} = 0$.
Therefore $u_k(\zeta^k-\zeta^{-k})=0$ for all $k\ge0$. Taking $\zeta$ to be a root of unity of precise order~$2e$, it follows that $u_k=0$ for all $k\ge0$ which are not multiples of~$e$. This means that $u(z)=v(z^e)$ for a second analytic function~$v$. We conclude that locally near~$\tilde c$, both $\mu\circ\pi_1$ and $\mu\circ\pi_2$, and hence both $\pi_1$ and~$\pi_2$, factor through the map $z\mapsto z^e$. But since the map $\pi\colon\tilde C\to C$ is an isomorphism outside finitely many points, this is only possible with $e=1$. This proves (d).

Finally, (c) and (d) together imply (e).
\end{Proof}

\medskip
Now we turn to the proof of Theorem \ref{NS3}. To ease notation we rename the variables $(X_1,X_2)$ to $(X,Y)$. Let $A$ denote the ring of all power series in $\BC[[X,Y]]$ satisfying
\UseTheoremCounterForNextEquation
\begin{equation}\label{NormDefb}
\textstyle \bigl\Vert\sum_{j,k} a_{jk} X^jY^k\bigr\Vert \ :=\ \sum_{j,k} |a_{jk}| \ <\ \infty
\end{equation}
Equivalently this is the ring of all power series that converge absolutely on~$\BB^2$, or again the Banach algebra completion of $\BC[X,Y]$ with respect to $\Vert\ \Vert$.
For any real number $0<r<1$ we set 
\UseTheoremCounterForNextEquation
\begin{equation}\label{frDef}
f_r(X,Y)\ :=\ f(X,rY).
\end{equation}
Since by assumption $f$ vanishes nowhere on $\BB\times\BB^\circ$, the polynomial $f_r$ vanishes nowhere on~$\BB^2$. Thus it is non-zero on a whole neighborhood of $\BB^2$, and therefore invertible in~$A$. The main point in the proof will be Proposition \ref{NS3aMainProp}, which asserts that the functions $\frac{f}{f_r}$ for $r\nearrow1$ constitute something like an \emph{approximate identity} for the ideal $(f)$ (compare Mortini-von Renteln \cite[page 223]{MortiniVonRenteln1989}). 
Abbreviate $\partial_Y := \frac{\partial}{\partial Y}$. For any $\ell\ge0$ write
\UseTheoremCounterForNextEquation
\begin{equation}\label{gDef}
g_\ell\ :=\ \frac{\partial_Y^\ell f}{f}
\ =\ \sum_{k\ge0} g_{\ell k}(X) Y^k
\ =\ \sum_{k,j\ge0} a_{\ell kj} X^j Y^k.
\end{equation}
By Proposition \ref{PlaneDClass} (c) the polynomial $f$ has degree $d_2$ with respect to~$Y$; hence $g_\ell=0$ for all $\ell>d_2$.

\begin{Lem}\label{NS3aLem1}
For all $\ell$ and $k$ and all $x\in\partial\BB$ we have 
$$|g_{\ell k}(x)|\ =\ O\bigl((k+1)^{\ell-1}\bigr),$$
where the implicit constant is independent of $\ell$, $k$, and~$x$.
\end{Lem}

\begin{Proof}
By assumption the polynomial $f(X,0)$ vanishes nowhere on~$\BB$. For any $x\in\partial\BB$ Proposition \ref{PlaneDClass} (e) therefore shows that
$$f(x,Y)\ =\ f(x,0)\cdot\!\!\!\prod_{\tilde c\in\pi_1^{-1}(x)}\Bigl(1-\frac{Y}{\pi_2(\tilde c)}\Bigr).$$
By the Leibniz formula it follows that 
$$\partial_Y^\ell f(x,Y)\ =\ f(x,0)\cdot\!\!\!
\sum_{{\scriptstyle I\subset\pi_1^{-1}(x)\atop \scriptstyle |I|=\ell}}
\,\prod_{\tilde c\in I}\ \Bigl(-\frac{1}{\pi_2(\tilde c)}\Bigr) \cdot\!\!
\prod_{\tilde c\in\pi_1^{-1}(x)\setminus I} \!\!\Bigl(1-\frac{Y}{\pi_2(\tilde c)}\Bigr).$$
Thus
\begin{eqnarray*}
g_\ell(x,Y) &\!\!=\!\!& 
\sum_{{\scriptstyle I\subset\pi_1^{-1}(x)\atop\scriptstyle  |I|=\ell}}
\,\prod_{\tilde c\in I}\ 
\biggl(\Bigl(-\frac{1}{\pi_2(\tilde c)}\Bigr)\!\!\Bigm/\!\!
\Bigl(1-\frac{Y}{\pi_2(\tilde c)}\Bigr) \biggr) \\
&\!\!=\!\!& 
\sum_{{\scriptstyle I\subset\pi_1^{-1}(x)\atop\scriptstyle  |I|=\ell}}
\prod_{\tilde c\in I}\ 
\biggl(\sum_{j\ge1} \frac{-Y^{j-1}}{\pi_2(\tilde c)^j}\biggr) \\
&\!\!=\!\!& 
\sum_{{\scriptstyle I\subset\pi_1^{-1}(x)\atop\scriptstyle |I|=\ell}}
\;\sum_{j\colon I\to\BZ^{\ge1}}\,
\!\!(-1)^\ell\cdot \prod_{\tilde c\in I}\ 
\frac{Y^{j(\tilde c)-1}}{\pi_2(\tilde c)^{j(\tilde c)}} \\
&\!\!=\!\!& 
\sum_{{\scriptstyle j\colon\pi_1^{-1}(x)\to\BZ^{\ge0}
\atop\scriptstyle  |\{\tilde c\,|\,j(\tilde c)>0\}|=\ell}}
\!\!\!\!(-1)^\ell\cdot\biggl(\,\prod_{\tilde c\in\pi_1^{-1}(x)}
\frac{1}{\pi_2(\tilde c)^{j(\tilde c)}} \biggr)
\cdot Y^{\sum_{\tilde c}j(\tilde c)-\ell} \\
\end{eqnarray*}
For any $k\ge0$ we therefore have 
\UseTheoremCounterForNextEquation
\begin{equation}\label{NS3aLem1a}
g_{\ell k}(x) \ =\ 
\sum_{j\in N_{\ell,k}(x)}
\!\!(-1)^\ell\cdot\!\!\!\!\prod_{\tilde c\in\pi_1^{-1}(x)}
\frac{1}{\pi_2(\tilde c)^{j(\tilde c)}}
\end{equation}
where $N_{\ell,k}(x)$ denotes the set of all maps $j\colon\pi_1^{-1}(x)\to\BZ^{\ge0}$ satisfying $|\{\tilde c\,|\,j(\tilde c)>0\}|=\ell$ and $\sum_{\tilde c}j(\tilde c)-\ell=k$.
But by Proposition \ref{PlaneDClass} (b), for all $\tilde c\in \pi_1^{-1}(x)$ we have $|\pi_2(\tilde c))|=1$. Since $|\pi_1^{-1}(x)|=d_2$ by Proposition \ref{PlaneDClass} (e), it follows that  
$$|g_{\ell k}(x)|\ \le\ |N_{\ell,k}(x)|
\ =\ \binom{d_2}{\ell}\cdot\binom{\ell+k-1}{\ell-1}.$$
Here the right hand side is $0$ for $\ell>d_2$, and a polynomial of degree $\ell-1$ in $k$ otherwise. The desired estimate follows.
\end{Proof}

\begin{Lem}\label{NS3aLem2}
There exists $M>0$ such that for all $\ell$ and $k$ we have 
$$\sum_{j>kM} |a_{\ell kj}| \ =\ O\bigl((k+1)^{\ell-1}\bigr),$$
where the implicit constant is independent of $\ell$ and $k$.
\end{Lem}

\begin{Proof}
Recall that by Proposition \ref{PlaneDClass} (b) and (d), for any $x\in\partial\BB$ and any $\tilde c\in\pi_1^{-1}(x)$ we have $|\pi_2(\tilde c)|=1$, and $\pi_1$ is unramified at~$\tilde c$. By continuity it follows that for any $x$ in a suitable neighborhood of $\partial\BB$ and any $\tilde c\in\pi_1^{-1}(x)$ we have $|\pi_2(\tilde c)|\ge\frac{1}{2}$, and $\pi_1$ is unramified at~$\tilde c$. Also, since $f(x,0)$ vanishes nowhere on~$\partial\BB$, the same holds in a neighborhood. Choose $\rho>1$ such that the slightly larger circle $\rho\cdot\partial\BB$ is contained in both neighborhoods. Then for all $x\in\rho\cdot\partial\BB$, the formula (\ref{NS3aLem1a}) remains true. Using $|\pi_2(\tilde c)|\ge\frac{1}{2}$ it now yields the estimate
$$|g_{\ell k}(x)|\ \le\ |N_{\ell,k}(x)|\cdot 2^{\ell+k}
\ =\ O\bigl((k+1)^{\ell-1} 2^k\bigr).$$
Plugging this into the Cauchy integral formula, for any $j\ge0$ we deduce that
$$|a_{\ell kj}| \ =\ 
\biggl| \frac{1}{2\pi i} \int_{\rho\cdot\partial\BB} g_{\ell k}(x)\, \frac{dx}{x^{j+1}} \biggr| 
\ \le \ 
\frac{1}{2\pi} \int_{\rho\cdot\partial\BB} |g_{\ell k}(x)|\, 
\frac{|dx|}{\rho^{j+1}}
\ = \ O\bigl((k+1)^{\ell-1} 2^k\rho^{-j}\bigr).$$
Summing the geometric series $\sum_{j>kM}\rho^{-j} = \rho^{-kM}/(\rho-1)$, we obtain the estimate
$$\sum_{j>kM} |a_{\ell kj}| \ =\ O\bigl((k+1)^{\ell-1}2^k\rho^{-kM}\bigr).$$
Any $M>0$ with $\rho^M\ge2$ thus has the desired property.
\end{Proof}

\begin{Lem}\label{NS3aLem3}
For all $\ell$ and $k$ we have 
$$\Vert g_{\ell k}\Vert\ =\ O\bigl((k+1)^{\ell-\frac{1}{2}}\bigr),$$
where the implicit constant is independent of $\ell$ and $k$.
\end{Lem}

\begin{Proof}
By (\ref{gDef}) we have $g_{\ell k} = \sum_{j\ge0} a_{\ell kj} X^j$.
Thus by the Parseval identity $(\sum_{j\ge0} |a_{\ell kj}|^2)^{\frac{1}{2}}$ is equal to the $L^2$-norm of the function $g_{\ell k}|\partial\BB$. It is therefore less than or equal to the $L^\infty$-norm of $g_{\ell k}|\partial\BB$. By Lemma \ref{NS3aLem1} it follows that 
$$\biggl(\sum_{j\ge0} |a_{\ell kj}|^2\biggr)^{\frac{1}{2}} 
\ =\ O\bigl((k+1)^{\ell-1}\bigr).$$
For any fixed $M>0$, with the Cauchy-Schwarz inequality we deduce that 
$$\sum_{0\le j\le kM} |a_{\ell kj}| \ \le\ 
\biggl(\,\sum_{0\le j\le kM}\!\! |a_{\ell kj}|^2\biggr)^{\frac{1}{2}}
\cdot \biggl(\,\sum_{0\le j\le kM}\!\! 1\biggr)^{\frac{1}{2}}
\ \le\ O\bigl((k+1)^{\ell-1}\bigr)\cdot (kM+1)^{\frac{1}{2}}.$$
Combining this with the estimate from Lemma \ref{NS3aLem2} for a suitable choice of $M>0$ we conclude that
$$\Vert g_{\ell k}\Vert\ =\ \sum_{j\ge0} |a_{\ell kj}|
\ =\ O\bigl((k+1)^{\ell-1}\bigr)\cdot (kM+1)^{\frac{1}{2}}
\ =\ O\bigl((k+1)^{\ell-\frac{1}{2}}\bigr),$$
as desired.
\end{Proof}

\begin{Lem}\label{NS3aLem4}
For all $\ell$ and all $0<r<1$ we have 
$$\Vert g_\ell(X,rY)\Vert\ =\ O\bigl((1-r)^{-\ell-\frac{1}{2}}\bigr),$$
where the implicit constant is independent of $\ell$ and $r$.
\end{Lem}

\begin{Proof}
By (\ref{gDef}) and Lemma \ref{NS3aLem3} we have 
$$\Vert g_\ell(X,rY)\Vert
\ =\ \Bigl\Vert\sum_{k\ge0} g_{\ell k}(X) r^kY^k\Bigr\Vert
\ =\ \sum_{k\ge0} \Vert g_{\ell k}(X)\Vert r^k
\ =\ O\Bigl(\sum_{k\ge0} (k+1)^{\ell-\frac{1}{2}}r^k\Bigr).$$
The usual integral estimate and the substitution $t=\frac{-s}{\log r}$ yield
\begin{eqnarray*}
\sum_{k\ge0} (k+1)^{\ell-\frac{1}{2}}r^k
&\!\!\le\!\!& \int_0^\infty (t+1)^{\ell-\frac{1}{2}}r^{t-1}\,dt \\
&\!\!=\!\!& \int_0^\infty \Bigl(1-\frac{s}{\log r}\Bigr)^{\ell-\frac{1}{2}}
\cdot\frac{e^{-s}}{r}\cdot\frac{-ds}{\log r} \\
&\!\!=\!\!& \Bigl(\frac{-1}{\log r}\Bigr)^{\ell+\frac{1}{2}}\cdot\frac{1}{r} \cdot\int_0^\infty (s-\log r)^{\ell-\frac{1}{2}}\cdot e^{-s}\,ds.
\end{eqnarray*}
Since $\frac{-1}{\log r}\sim \frac{1}{1-r}$ for $r\nearrow1$, and the integral converges to a finite value, the desired estimate follows.
\end{Proof}

\begin{Lem}\label{NS3aLem5}
For $r\nearrow1$ we have 
$$\Bigl\Vert\frac{f}{f_r}\Bigr\Vert\ =\ O\bigl((1-r)^{-\frac{1}{2}}\bigr).$$
\end{Lem}

\begin{Proof}
The formal Taylor expansion of $f(X,Y)$ at $Y=Y_0$ reads
$$f(X,Y)\ =\ \sum_{0\le\ell\le d_2} 
(\partial_Y^\ell f)(X,Y_0)\cdot\frac{(Y-Y_0)^\ell}{\ell!}.$$
Substituting $Y_0=rY$ and using (\ref{gDef}) we deduce that 
\begin{eqnarray*}
\frac{f}{f_r}(X,Y)
&\!\!=\!\!& \sum_{0\le\ell\le d_2} 
\frac{(\partial_Y^\ell f)(X,rY)}{f(X,rY)}\cdot\frac{(Y-rY)^\ell}{\ell!} \\
&\!\!=\!\!& \sum_{0\le\ell\le d_2} 
g_\ell(X,rY)\cdot\frac{(1-r)^\ell}{\ell!}\cdot Y^\ell.
\end{eqnarray*}
Since $\bigl\Vert g_\ell(X,rY)\cdot(1-r)^\ell\bigr\Vert = O\bigl((1-r)^{-\frac{1}{2}}\bigr)$ by Lemma \ref{NS3aLem4}, the estimate follows.
\end{Proof}

\begin{Prop}\label{NS3aMainProp}
For $r\nearrow1$ we have 
$$\Bigl\Vert \frac{f^2}{f_r}-f\Bigr\Vert\ =\ O\bigl((1-r)^{\frac{1}{2}}\bigr).$$
\end{Prop}

\begin{Proof}
Combining $\Vert f-f_r\Vert=O(1-r)$ and Lemma \ref{NS3aLem5} yields
$$\Bigl\Vert \frac{f^2}{f_r}-f\Bigr\Vert\ =\ 
\Bigl\Vert (f-f_r)\frac{f}{f_r}\Bigr\Vert\ \le\ 
\Vert f-f_r\Vert\cdot\Bigl\Vert\frac{f}{f_r}\Bigr\Vert\ =\ 
O(1-r)\cdot O\bigl((1-r)^{-\frac{1}{2}}\bigr),$$
as desired.
\end{Proof}

\begin{Thm}\label{NS3a}
{\bf(=\ \ref{NS3})}\ 
In the case (d) of Proposition \ref{PlaneClass}, for all $k\ge1$ we have 
$$\overline{(f^k)} \ =\ \widetilde{(f^k)} \ =\ (f).$$
\end{Thm}

\begin{Proof}
For each $0<r<1$, since $\frac{1}{f_r}\in A$, there exists a polynomial $h_r\in\BC[X,Y]$ with $\bigl\Vert h_r-\frac{1}{f_r}\bigr\Vert\le (1-r)^{\frac{1}{2}}$. Using Proposition \ref{NS3aMainProp} we deduce that
$$\bigl\Vert f^2h_r-f\bigr\Vert
\ \le\ \Vert f^2\Vert\cdot\Bigl\Vert h_r-\frac{1}{f_r}\Bigr\Vert
+ \Bigl\Vert\frac{f^2}{f_r}-f\Bigr\Vert
\ =\ O\bigl((1-r)^{\frac{1}{2}}\bigr).$$
As this tends to $0$ for $r\nearrow1$, it follows that $f$ lies in the closure of the ideal $(f^2)\subset\BC[X,Y]$. 
By induction on $k$, as in the proof of Lemma \ref{1var3}, we deduce that $f\in\overline{(f^k)}$ for all $k\ge1$. Using Proposition \ref{ITildeProp} (c) and (b) we conclude that
$$(f)\ \subset\ \overline{(f^k)} \ \subset\ \widetilde{(f^k)} 
\ \subset\ I(M((f^k)))\ =\ (f).$$
The desired equalities follow.
\end{Proof}


\section{General case in the plane}
\label{Dim2}

We keep $n=2$.

\begin{Lem}\label{Dim2Lem1}
Consider any non-zero polynomial $f\in\BC[\UX]$, any point $\Uz=(z_1,z_2)\in\BB^2$, and any closed $\Fm_\Uz$-primary ideal $J\subsetneqq\BC[\UX]$. Then
$$(f)\cap \bigcap_{r\ge1} \bigl(fJ+\overline{\Fm_\Uz^r}\,\bigr)
\ =\ \left\{\begin{array}{cl}
(f) & \hbox{if $\Uz\in\partial\BB\times\BB^\circ$ and $f\in(X_1-z_1)$,} \\[3pt]
(f)& \hbox{if $\Uz\in\BB^\circ\times\partial\BB$ and $f\in(X_2-z_2)$,} \\[3pt]
(f)& \hbox{if $\Uz\in(\partial\BB)^2$ and $f\in\Fm_\Uz$,} \\[3pt]
fJ & \hbox{otherwise.}
\end{array}\right.$$
\end{Lem}

\begin{Proof}
Call the left hand side $I$ and note that we always have $I\supset fJ$.

Suppose first that $\Uz\in(\BB^\circ)^2$. Then for each $r\ge1$ we have $\overline{\Fm_\Uz^r}=\Fm_\Uz^r$ by Theorem \ref{NS2}. Using the primary decomposition we can write $fJ=(f)\cap J'$ for some $\Fm_\Uz$-primary ideal~$J'$. Then for some $r\ge1$ we have $\Fm_\Uz^r\subset J'$ and hence $fJ+\overline{\Fm_\Uz^r} = ((f)\cap J')+\Fm_\Uz^r \subset J'$. Thus $(f)\cap\bigl(fJ+\overline{\Fm_\Uz^r}\bigr) \subset (f)\cap J' = fJ$. This implies that $I\subset fJ$ and hence $I=fJ$, as desired.

Suppose next that $\Uz\in\partial\BB\times\BB^\circ$. Then for each $r\ge1$ we have $\overline{\Fm_\Uz^r}=\bigl(X_1-z_1,(X_2-z_2)^r\bigr)$ by Theorem \ref{NS2}. Since $J$ is $\Fm_\Uz$-primary and closed, we also have $J=\bigl(X_1-z_1,(X_2-z_2)^s\bigr)$ for some $s\ge1$. For all $g\in\BC[\UX]$ we deduce that
\begin{eqnarray*}
fg \in I &\Longleftrightarrow& \forall r\ge1\colon \ fg\ \in\ 
f\cdot\bigl(X_1-z_1,(X_2-z_2)^s\bigr)+\bigl(X_1-z_1,(X_2-z_2)^r\bigr) \\
&\Longleftrightarrow& \forall r\ge1\colon \ (fg)(z_1,X_2)\ \in\ \bigl(f(z_1,X_2)\cdot(X_2-z_2)^s,(X_2-z_2)^r\bigr) \ \subset\ \BC[X_2].
\end{eqnarray*}
If $f(z_1,X_2)=0$, or equivalently $f\in(X_1-z_1)$, this condition always holds; hence in this case $I=(f)$. Otherwise, taking $r-s$ greater than the maximal power of $(X_2-z_2)$ dividing $f(z_1,X_2)$, we find that
\begin{eqnarray*}
fg \in I &\Longleftrightarrow& 
g(z_1,X_2)\ \in\ \bigl((X_2-z_2)^s\bigr) \ \subset\ \BC[X_2] \\
&\Longleftrightarrow& 
g\ \in\ \bigl(X_1-z_1,(X_2-z_2)^s\bigr) \ =\ J.
\end{eqnarray*}
In this case we therefore have $I=fJ$.

The case $\Uz\in\BB^\circ\times\partial\BB$ follows by symmetry from the preceding case.

Suppose finally that $\Uz\in(\partial\BB)^2$. Then for each $r\ge1$ we have $\overline{\Fm_\Uz^r}=\Fm_\Uz$ by Theorem \ref{NS2}. Since $J$ is $\Fm_\Uz$-primary and proper and closed, we also have $J=\Fm_\Uz$. 
We deduce that $fJ+\overline{\Fm_\Uz^r} = f\Fm_\Uz+\Fm_\Uz = \Fm_\Uz$, and hence $I=(f)\cap\Fm_\Uz$. It follows that $I=(f)$ if $f\in\Fm_\Uz$, and $I=f\Fm_\Uz=fJ$ otherwise.
This finishes the proof in all cases.
\end{Proof}

\begin{Lem}\label{Dim2Lem2}
In all cases of Lemma \ref{Dim2Lem1} we have 
$$(f)\cap \bigcap_{r\ge1} \bigl(fJ+\overline{\Fm_\Uz^r}\,\bigr)
\ =\ (f)\cap \overline{(fJ)}.$$
\end{Lem}

\begin{Proof}
The ideals $fJ+\overline{\Fm_\Uz^r}$ are all closed by Lemma \ref{FiniteCodimLem}, hence they contain $\overline{(fJ)}$. Calling the left hand side $I$, we thus always have $I\supset (f)\cap \overline{(fJ)}\supset fJ$. We are therefore done if $I=fJ$. It remains to consider the cases of Lemma \ref{Dim2Lem1} where $I=(f)$.

Suppose first that $\Uz\in\partial\BB\times\BB^\circ$ and $f\in(X_1-z_1)$. 
Write $f=(X_1-z_1)f_1$ with $f_1\in\BC[\UX]$. Since $(X_1-z_1)^r\in J$ for some $r\ge1$, we have $(X_1-z_1)^{r+1}f_1 \in fJ$. But by Lemma \ref{1var3} and the isometric embedding $\BC[X_1]\into\BC[\UX]$ we have $(X_1-z_1)\in \overline{((X_1-z_1)^{r+1})}$. Therefore $f=(X_1-z_1)f_1 \in \overline{((X_1-z_1)^{r+1})f_1} \subset \overline{(fJ)}$ and hence $(f)\cap \overline{(fJ)}=(f)=I$, as desired.

The case $\Uz\in\BB^\circ\times\partial\BB$ and $f\in(X_2-z_2)$ follows by symmetry from the preceding case.

Suppose lastly that $\Uz\in(\partial\BB)^2$ and $f\in\Fm_\Uz$. Recall from the proof of Lemma \ref{Dim2Lem1} that $J=\Fm_\Uz$ in this case.
Write $f=(X_1-z_1)f_1+(X_2-z_2)f_2$ with $f_1,f_2\in\BC[\UX]$. By the same argument as in the proof of Lemma \ref{1var3}, for each $i=1,2$ and all $m\ge1$ we have
$$1\ \equiv\ \frac{1}{mz_i^{m-1}}\cdot\frac{X_i^m-z_i^m}{X_i-z_i}
\ =\ \frac{X_i^{m-1}+X_i^{m-2}z_i+\ldots+z_i^{m-1}}{mz_i^{m-1}}
\ \ \hbox{modulo}\ \ (X_i-z_i).$$
Modulo $f\Fm_\Uz=fJ$ we therefore have
\begin{eqnarray*}
f &\!\!\equiv\!\!&
\frac{1}{mz_1^{m-1}}\cdot\frac{X_1^m-z_1^m}{X_1-z_1}\cdot
\frac{1}{mz_2^{m-1}}\cdot\frac{X_2^m-z_2^m}{X_2-z_2}\cdot f \\
&\!\!=\!\!&
\frac{X_1^m-z_1^m}{mz_1^{m-1}}\cdot
\frac{1}{mz_2^{m-1}}\cdot\frac{X_2^m-z_2^m}{X_2-z_2}\cdot f_1
+ \frac{1}{mz_1^{m-1}}\cdot\frac{X_1^m-z_1^m}{X_1-z_1}\cdot
\frac{X_2^m-z_2^m}{mz_2^{m-1}}\cdot f_2.
\end{eqnarray*}
Here, since $|z_i|=1$, we have $\bigl\Vert\frac{1}{mz_i^{m-1}}(X_i^m-z_i^m)\bigr\Vert = \frac{1}{m}(1+1)\to0$ for $m\to\infty$, whereas $\bigl\Vert\frac{1}{mz_i^{m-1}}\cdot\frac{X_i^m-z_i^m}{X_i-z_i}\bigr\Vert = \frac{1}{m}\bigl\Vert X_i^{m-1}+X_i^{m-2}z_i+\ldots+z_i^{m-1}\bigr\Vert 
= 1$ for all~$m$. In the limit we thus deduce that $f\equiv0$ modulo $\overline{(fJ)}$. Therefore $(f)\cap \overline{(fJ)}=(f)=I$, as desired.
\end{Proof}

\begin{Lem}\label{Dim2Lem3}
In all cases of Lemma \ref{Dim2Lem1} we have 
$$(f)\cap\widetilde{(fJ)}\ =\ (f)\cap \overline{(fJ)}.$$
\end{Lem}

\begin{Proof}
Since $J$ is $\Fm_\Uz$-primary, for all points $\Uz'\in\BB^2$ distinct from $\Uz$ and all $r\ge1$ we have 
$(f)\cap \bigl(fJ+\overline{\Fm_{\smash{\Uz'}}^r}\,\bigr)=(f)$. Thus the ideal on the left hand side of Lemma \ref{Dim2Lem2} is just $(f)\cap\widetilde{(fJ)}$.
\end{Proof}

\begin{Thm}\label{NS4b}
For any ideal $I\subset\BC[X_1,X_2]$ we have $\OI=\TI$.
\end{Thm}

\begin{Proof}
If $I=0$, we have $M(I)=\BB^2$ and hence $I(M(I))=0$. With Proposition \ref{ITildeProp} (b) we deduce that $\TI=0=\OI$, and the theorem follows. Henceforth we assume that $I\not=0$. 

By Proposition \ref{ITildeProp} (c) we can replace $I$ by $\OI$ without changing~$\TI$. Thus without loss of generality we assume that $I$ is closed.

Write $I=fJ$ for a non-zero polynomial $f\in\BC[\UX]$ and an ideal of finite codimension $J\subset\BC[\UX]$. Write $f=\prod_{i=1}^m f_i^{k_i}$ with mutually non-associate irreducible polynomials $f_i\in\BC[\UX]$ and exponents $k_i\ge1$. 

\begin{Claim}\label{NS4bLem1}
For all $i$ we have $(f_i^{k_i})=\overline{(f_i^{k_i})}=\widetilde{(f_i^{k_i})}$.
\end{Claim}

\begin{Proof}
By Theorem \ref{NS1} we have $I(M(I)) = \Rad(I)$. For each $i$ the set  $M(I)$ therefore contains a Zariski dense subset of the curve defined by $f_i$. Thus $M((f_i))$ is infinite, and so $f_i$ does not belong to the case (e) of Proposition \ref{PlaneClass}. Set $k_i':=k_i$ if $f_i$ belongs to the case \ref{PlaneClass} (a), and $k_i':=1$ if it belongs to the cases \ref{PlaneClass} (b)--(d). Then by Propositions \ref{PlaneA} and \ref{PlaneBC} and Theorem \ref{NS3a} we have 
$$\overline{(f_i^{k_i})} \ =\ \widetilde{(f_i^{k_i})} \ =\ (f_i^{k_i'}).$$
By the continuity of addition and multiplication in $\BC[\UX]$ we deduce that 
$$\Bigl(\prod_{i=1}^m f_i^{k_i'}\Bigr) J \ =\ 
\Bigl(\prod_{i=1}^m \overline{(f_i^{k_i})}\Bigr) J \ \subset\  
\overline{\Bigl(\Bigl(\prod_{i=1}^m (f_i^{k_i})\Bigr) J\Bigr)} \ =\  
\OI \ =\  I\ =\ 
\Bigl(\prod_{i=1}^m f_i^{k_i}\Bigr) J.$$
Thus for all $i$ we have $k_i'\ge k_i$ and hence $k_i'=k_i$, and the claim follows.
\end{Proof}

\begin{Claim}\label{NS4bLem2}
We have $(f)=\overline{(f)}=\widetilde{(f)}$.
\end{Claim}

\begin{Proof}
For each $i$ we have $(f)\subset(f_i^{k_i})$. Using Proposition \ref{ITildeProp} (c) and (d) and Claim \ref{NS4bLem1} we deduce that
$$(f)\ \subset\ \overline{(f)}\ \subset\ \widetilde{(f)}\ \subset\ 
\bigcap_{i=1}^m\widetilde{(f_i^{k_i})} \ =\ 
\bigcap_{i=1}^m(f_i^{k_i}) \ =\ (f),$$
from which the equalities follow.
\end{Proof}

\medskip
Now write $J=\bigcap_{\nu=1}^\ell J_\nu$ where the ideals $J_\nu\subsetneqq\BC[\UX]$ are $\Fm_{\Uz_\nu}$-primary for distinct points $\Uz_\nu\in\BC^2$. Then we also have 
$J=\prod_{\nu=1}^\ell J_\nu$. 

\begin{Claim}\label{NS4bLem3}
For each $\nu$ the ideal $fJ_\nu$ is closed.
\end{Claim}

\begin{Proof}
Since $(f)$ is closed by Claim \ref{NS4bLem2}, we have $\overline{(fJ_\nu)}=fJ_\nu^+$ for some ideal $J_\nu^+$ containing~$J_\nu$. Abbreviate $J' := \prod_{\nu'\not=\nu} J_{\nu'}$, so that $fJ_\nu J'=fJ=I$. 
Using continuity of addition and multiplication, we calculate
$$fJ_\nu^+J' \ =\ \overline{(fJ_\nu)}J' \ \subset\
\overline{(fJ_\nu)J'} \ =\ \OI \ =\  I\ =\ fJ_\nu J'.$$
Dividing by $f$ and using the fact that $J_\nu$ and $J'$ have disjoint support, we deduce that $J_\nu^+\subset J_\nu$. Thus $J_\nu^+=J_\nu$ and hence $\overline{(fJ_\nu)}=fJ_\nu^+=fJ_\nu$, as desired.
\end{Proof}

\begin{Claim}\label{NS4bLem4}
For each $\nu$ the ideal $J_\nu$ is closed.
\end{Claim}

\begin{Proof}
By continuity and Claim \ref{NS4bLem3} we have $fJ_\nu\subset f\cdot\overline{J_\nu} \subset \overline{(fJ_\nu)} = fJ_\nu$. Dividing by $f$ yields $J_\nu\subset\overline{J_\nu}\subset J_\nu$.
\end{Proof}

\begin{Claim}\label{NS4bLem5}
For each $\nu$ we have $\widetilde{(fJ_\nu)}= fJ_\nu$.
\end{Claim}

\begin{Proof}
Since $J_\nu$ is closed, Theorem \ref{NS2} implies that the associated point $\Uz_\nu$ lies in~$\BB^2$. We can therefore apply Lemma \ref{Dim2Lem3} to $f$ and~$J_\nu$, yielding
$$(f)\cap\widetilde{(fJ_\nu)} \ =\  (f)\cap \overline{(fJ_\nu)}.$$
By Proposition \ref{ITildeProp} (d) and Claim \ref{NS4bLem2}
we also have $\widetilde{(fJ_\nu)} \subset \widetilde{(f)} = (f)$; hence the left hand side is just $\widetilde{(fJ_\nu)}$. On the other hand, since $fJ_\nu$ is closed by Claim \ref{NS4bLem3}, the right hand side is just $fJ_\nu$. Thus the desired equality follows.
\end{Proof}

\medskip
\noindent{\bf End of Proof.} 
Finally, using Proposition \ref{ITildeProp} (b) and (d) and Claim \ref{NS4bLem5} we conclude that
$$fJ\ \subset\ \widetilde{(fJ)}\ \subset\ 
\bigcap_{\nu=1}^\ell \widetilde{(fJ_\nu)}
\ =\ \bigcap_{\nu=1}^\ell fJ_\nu \ =\ fJ.$$
Thus $\widetilde{(fJ)}=fJ$. Since $fJ=I=\OI$, it follows that $\OI=\TI$.
This finishes the proof of Theorem \ref{NS4b}.
\end{Proof}


\section{Examples}
\label{Ex}

\begin{Ex}\label{Ex1}
\rm For $f:=X_1+X_2-2\in\BC[X_1,X_2]$ we have $M((f))=\{(1,1)\}$. This is therefore the case (e) of Proposition \ref{PlaneClass} with non-empty finite $M((f))\subset(\partial\BB)^2$. In this case Theorem \ref{NS4} asserts that $\overline{(f)} = (X_1-1,X_2-1)$. 
We give a direct proof of this fact.
\end{Ex}

\begin{Proof}
Since $X_1+X_2\equiv 2$ modulo $(f)$, for every $m\ge0$ we have 
$$g_m \ :=\ (X_1+X_2)^{2m}(X_1-X_2) \ \equiv\ 2^m(X_1-X_2)\ \ \hbox{modulo}\ \ (f).$$
Expanding $g_m$ with the binomial formula yields
\begin{eqnarray*}
g_m &\!\!=\!\!& \sum_{i=0}^{2m} \binom{2m}{i} X_1^i X_2^{2m-i} (X_1-X_2) \\
&\!\!=\!\!& X_1^{2m+1} + \sum_{i=1}^{2m} \left(\binom{2m}{i-1}-\binom{2m}{i}\right) X_1^i X_2^{2m+1-i} - X_2^{2m+1}.
\end{eqnarray*}
Using a telescoping sum, we can therefore determine its norm as
\begin{eqnarray*}
\Vert g_m\Vert &\!\!=\!\!& 1 + \sum_{i=1}^{2m} \left|\binom{2m}{i-1}-\binom{2m}{i}\right| + 1 \\
&\!\!=\!\!& 2 + 2\sum_{i=1}^{m} \left(\binom{2m}{i}-\binom{2m}{i-1}\right) \ =\ 2\binom{2m}{m}.
\end{eqnarray*}
It is a well-known consequence of Stirling's formula that $\binom{2m}{m} = O(2^m/\sqrt{m})$ as $m\to\infty$. Thus $2^{-m}g_m\to0$ for $m\to\infty$. 
Since $X_1-X_2\equiv2^{-m}g_m$ modulo $(f)$, it follows that $X_1-X_2\in\overline{(f)}$. Therefore $X_1-1 = \frac{1}{2}((X_1-X_2)+f)\in\overline{(f)}$, and similarly $X_2-1\in\overline{(f)}$. Thus $(X_1-1,X_2-1)\subset\overline{(f)}$. The reverse inclusion follows from the fact that $\overline{(f)}\subset I(M((f)))$.
\end{Proof}

\begin{Ex}\label{Ex2}
\rm Fix a real number $0<w<1$ and consider the polynomial
$$f\ :=\ 1+wX+wY+XY \ \in \ \BC[X,Y].$$
Let $C$ denote the closure of the affine curve $V(f)$ in~$\hat\BC^2$. Its image under the M\"obius transformation from Section \ref{PlaneD} is the real hyperbola with the affine equation $UV=\frac{1+w}{1-w}$. Since $\frac{1+w}{1-w}>0$, for any complex point on this hyperbola with one coordinate in the upper half plane~$\BH$, the other coordinate lies in the lower half plane~$-\BH$. This implies that $M((f))$ is infinite and contained in $(\partial\BB)^2$, so that we are in the case (d) of Proposition \ref{PlaneClass}.

Among the functions $g_\ell$ from (\ref{gDef}), the only non-trivial one is now
$$g_1\ =\ \frac{w+X}{(1+wX)+(w+X)Y}
\ =\ \sum_{k\ge0} g_{1k}(X) Y^k$$ 
with 
$$g_{1k}\ =\ -\Bigl(-\frac{w+X}{1+wX}\Bigr)^{k+1}.$$
Expressing its coefficients with the Cauchy integral formula yields
$$a_{1kj} \ =\ 
\frac{-1}{2\pi} \int_{-\pi}^{\pi} \Bigl(-\frac{w+e^{it}}{1+we^{it}}\Bigr)^{k+1}\, \frac{dt}{e^{ijt}}.$$
Since $|w|<1$, we can write
$$-\frac{w+e^{it}}{1+we^{it}}\ =\ e^{i\phi(t)}$$
for a smooth function $\phi\colon\BR\to\BR$, and deduce that
$$a_{1kj} \ =\ \frac{-1}{2\pi} \int_{-\pi}^{\pi} 
e^{i(\phi(t)(k+1)-tj)}\,dt.$$
This oscillatory integral can be estimated with the method of stationary phase. It turns out that critical points exist if and only if
\UseTheoremCounterForNextEquation
\begin{equation}\label{interval}
\frac{1-w}{1+w}\ \le\ \frac{j}{k+1}\ \le\ \frac{1+w}{1-w}.
\end{equation}
Their contribution to the integral then comes out to be
$$\theta_{kj}\cdot \Bigl(\frac{\pi j}{2}\Bigr)^{-\frac{1}{2}}
\cdot \Bigl(\frac{j}{k+1} - \frac{1-w}{1+w}\Bigr)^{-\frac{1}{4}}
\cdot \Bigl(\frac{1+w}{1-w} - \frac{j}{k+1}\Bigr)^{-\frac{1}{4}}$$
with a wildly oscillating phase term $\theta_{kj}\in[-1,1]$.
Summing $|a_{1kj}|$ over $j$ in the range (\ref{interval}) shows that the order of magnitude of $\Vert g_{1k}\Vert$ is precisely $(k+1)^{\frac{1}{2}}$, and so Lemma \ref{NS3aLem3} is sharp in this case.
\end{Ex}

\begin{Rem}\label{Abund1}
\rm The case (d) of Proposition \ref{PlaneClass} is not particularly rare. For example, for any relatively prime integers $d_1$, $d_2\ge1$ the polynomial $X_1^{d_1} X_2^{d_2}-1$ satisfies the conditions and defines a smooth irreducible curve $C$ of bidegree $(d_1,d_2)$. Likewise, whenever $f$ satisfies the conditions, so does every irreducible factor of $f(X_1^{d_1},X_2^{d_2})$. 
Moreover, the conditions are invariant under applying a M\"obius transformation that maps $\BB$ to itself in each variable separately. 
Furthermore, the conditions are preserved under small deformations at least for smooth curves, by the following result.
\end{Rem}

\begin{Prop}\label{Abund2}
Let $S$ be an algebraic variety over $\BR$, and let $\CC\to S_\BC$ be an algebraic family of smooth irreducible curves in~$(\BP^1)^2$, such that the family obtained by applying $\mu\times\mu$ is defined over~$\BR$. 
Assume that in the fiber over some point $s_0\in S(\BR)$ the intersection $\CC_{s_0}(\BC)\cap\BB^2$ is infinite and contained in~$(\partial\BB)^2$. 
Then the same condition holds in the fiber over all points $s\in S(\BR)$ sufficiently near~$s_0$.
\end{Prop}

\begin{Proof}
Let $\CC'\to S$ be the family of curves in $(\BP^1)^2$ obtained by applying $\mu\times\mu$. Then Proposition \ref{PlaneDClass} for the curve at $s_0$ implies that the two projection maps
$$\xymatrix{\BP^1(\BR) & \CC'_{s_0}(\BR) \ar[l]_{\proj_1} \ar[r]^{\proj_2} & \BP^1(\BR)}$$
are unramified coverings, whose degrees are equal to the respective algebraic degrees $d_2$,~$d_1$. By smoothness the same is true over all nearby points $s\in S(\BR)$. 
Thus 
\UseTheoremCounterForNextEquation
\begin{equation}\label{Abund2a}
\CC'_s(\BC)\setminus\CC'_s(\BR) \ \subset\ (\BP^1(\BC)\setminus\BP^1(\BR))^2\ =\ (\BH\sqcup-\BH)^2.
\end{equation}
Topologically $\CC'_s(\BC)$ is a compact oriented surface of fixed genus and $\CC'_s(\BR)$ is the union of a fixed finite number of disjoint simple loops on it. By smoothness the families of both can be trivialized locally near~$s_0$.
The connected components of $\CC'_{s_0}(\BC)\setminus\CC'_{s_0}(\BR)$ are therefore in bijection with the connected components of $\CC'_s(\BC)\setminus\CC'_s(\BR)$ and deform continuously into the latter within $(\BH\sqcup-\BH)^2$. But the assumption $\CC_{s_0}(\BC)\cap(\BB^\circ)^2=\emptyset$ is equivalent to $\CC'_{s_0}(\BC)\cap\BH^2=\emptyset$. Thus none of the connected components of $\CC'_{s_0}(\BC)\setminus\CC'_{s_0}(\BR)$ meets~$\BH^2$, and so none of the connected components of $\CC'_s(\BC)\setminus\CC'_s(\BR)$ meets~$\BH^2$. 
Together with (\ref{Abund2a}) this implies that 
$$\CC'_s(\BC)\cap\bigl(\BH\sqcup\BP^1(\BR)\bigr)^2\ =\ \CC'_s(\BR).$$
Therefore $\CC_s(\BC)\cap\BB^2$ is infinite and contained in $(\partial\BB)^2$, as desired.
\end{Proof}



\end{document}